\documentstyle[amscd,amssymb,verbatim,12pt]{amsart}
\renewcommand{\mod}{\operatorname{mod}}

\newcommand{\und}{\underline}

\newcommand{\OO}{{\cal O}}

\newcommand{\coker}{\operatorname{coker}}
\newcommand{\DD}{{\cal D}}
\newcommand{\NN}{{\cal N}}

\newcommand{\hra}{\hookrightarrow}
\newcommand{\lan}{\langle}
\newcommand{\ran}{\rangle}

\newcommand{\CC}{{\cal C}}

\newcommand{\Spec}{\operatorname{Spec}}

\renewcommand{\P}{{\Bbb P}}

\newcommand{\Pic}{\operatorname{Pic}}

\newcommand{\eps}{\epsilon}

\renewcommand{\ker}{\operatorname{ker}}

\numberwithin{equation}{subsection}

\newtheorem{thm}{Theorem}[subsection]
\newtheorem{prop}[thm]{Proposition}
\newtheorem{lem}[thm]{Lemma}
\newtheorem{cor}[thm]{Corollary}
\newenvironment{rem}{\vspace{3mm}\noindent
{\bf Remark.}}{\vspace{3mm}}
\newenvironment{defi}{\vspace{3mm}\noindent
{\bf Definition.}}{\vspace{3mm}}

\newenvironment{ex}{\vspace{3mm}\noindent
{\bf Example.}}{\vspace{3mm}}

\newcommand{\Pf}{\noindent {\it Proof}}
\newcommand{\id}{\operatorname{id}}

\newcommand{\ov}{\overline}

\newcommand{\ra}{\rightarrow}
\newcommand{\HHom}{{\cal H}om}
\renewcommand{\AA}{{\cal A}}
\newcommand{\FF}{{\cal F}}
\newcommand{\TT}{{\cal T}}

\newcommand{\PP}{{\cal P}}

\newcommand{\QQ}{{\cal Q}}

\newcommand{\Hom}{\operatorname{Hom}}
\newcommand{\Autoeq}{\operatorname{Autoeq}}

\newcommand{\End}{\operatorname{End}}

\renewcommand{\a}{\alpha}
\renewcommand{\b}{\beta}
\newcommand{\om}{\omega}
\newcommand{\De}{\Delta}
\newcommand{\la}{\lambda}

\newcommand{\A}{{\Bbb A}}
\newcommand{\C}{{\Bbb C}}
\newcommand{\R}{{\Bbb R}}
\newcommand{\Z}{{\Bbb Z}}
\newcommand{\Q}{{\Bbb Q}}

\newcommand{\Ga}{\Gamma}

\newcommand{\wt}{\widetilde}

\newcommand{\sub}{\subset}
\newcommand{\ed}{\qed\vspace{3mm}}

\newcommand{\pra}{\operatorname{p-a}}
\newcommand{\apra}{\operatorname{a-p-a}}
\newcommand{\Du}{{\Bbb D}}

\title{Constant families of $t$-structures on derived categories of coherent sheaves}
\author{A. Polishchuk}
\address{Department of Mathematics, University of Oregon, Eugene, OR 97403}
\email{apolish@@uoregon.edu}
\thanks{Supported in part by the NSF grant DMS-0601034}

\begin{document}
\begin{abstract} 
We generalize the construction given in \cite{AP} of a ``constant" $t$-structure on the bounded derived category of coherent sheaves $D(X\times S)$ starting with a $t$-structure on $D(X)$. 
Namely, we remove smoothness and quasiprojectivity assumptions on $X$ and $S$
and work with $t$-structures that are not necessarily Noetherian but are close to Noetherian in the appropriate sense. The main new tool is the construction of induced $t$-structures
that uses unbounded derived categories of quasicoherent sheaves and relies on the results of 
\cite{AJS}. As an application of the ``constant" $t$-structures techniques we prove that every
bounded nondegenerate $t$-structure on $D(X)$ with Noetherian heart is invariant under
the action of a connected group of autoequivalences of $D(X)$. Also, we show that 
if $X$ is smooth then the only local $t$-structures on $D(X)$, i.e., those for which there exist
compatible $t$-structures on $D(U)$ for all open $U\sub X$, are the perverse $t$-structures
considered in \cite{Bezr}.
\end{abstract}
\maketitle

\bigskip

\centerline{\sc Introduction}

\medskip

Originally $t$-structures appeared in the context of derived categories of constructible sheaves
on a stratified space in the definition of perverse sheaves given in \cite{BBD}.
More recent studies led to interesting examples
of $t$-structures on bounded derived categories of coherent sheaves on algebraic varieties
(for example, in connection with the theory of stability conditions introduced by Bridgeland in \cite{Bridge1}; see also \cite{Bezr2} for examples relevant for representation theory). 
The present work is a continuation of \cite{AP} where we gave a construction of a $t$-structure on 
$D(X\times S)$, the bounded derived category of coherent sheaves on $X\times S$, 
starting with a $t$-structure on $D(X)$. This $t$-structure on $X\times S$ should be thought of as a constant family of $t$-structures over $S$ (we will often refer to it as a ``constant
$t$-structure"). Hopefully, it should serve as the first step
towards constructing nice moduli spaces for stable objects with respect to a stability condition on
$D(X)$ (see \cite{AP} for a discussion of this problem).

The main goal of this paper is to remove the smoothness assumption
that was imposed on $X$ and $S$ in \cite{AP}.
Moreover, we actually give an alternative construction even in the smooth case and remove
the assumption of boundedness with respect to the standard $t$-structure
in the results of \cite{AP}, sec.~2.7.
We still need the most nontrivial ingredient from \cite{AP} that gives the required $t$-structures
in the case $S=\P^r$. However, the remaining part of the construction is replaced by
a new method based on
the general procedure of ``inducing" a $t$-structure with respect
to a ``nice" functor (see Theorem \ref{induced-thm}). The geometric example of such a functor
relevant for the construction of constant $t$-structures
is the push-forward with respect to a finite morphism of finite Tor dimension.
The key idea is that it is much easier to construct $t$-structures in the unbounded derived categories
of quasicoherent sheaves $D_{qc}(X)$ 
because one can use arbitrary small coproducts. This idea was employed
effectively in Theorem A.1 of \cite{AJS} that shows that any pre-aisle stable
under all small coproducts and generated by a {\it set} of objects, extends to a $t$-structure (see section \ref{close-to-N-sec} for terminology). Of course, a random $t$-structure
on $D_{qc}(X)$ will not restrict to a $t$-structure on $D(X)$. However, if the two such categories
are related by a ``nice" functor $D_{qc}(X)\to D_{qc}(Y)$ then knowing that the $t$-structure on
$D_{qc}(Y)$ restricts to $D(Y)$ allows to deduce the same about the $t$-structure on $D_{qc}(X)$. 
Applying this approach we construct the constant $t$-structure on $D(X\times S)$ for
arbitrary $X$ and $S$ of finite type over a field (see Theorem \ref{const-t-str-thm}).

We also come up with several other improvements to \cite{AP}. First of all, in {\it loc. cit.}
we considered only Noetherian $t$-structures (i.e., $t$-structures with Noetherian heart).
In this paper we introduce {\it close to Noetherian} $t$-structures that are obtained
from Noetherian $t$-structures by tilting, and show that the construction of constant $t$-structures
goes through for them as well. A technical observation that facilitates such a generalization
is Theorem \ref{tilt-thm} stating that every pre-aisle, close to a Noetherian $t$-structure,
automatically extends to a $t$-structure. It is easy to see that in a reasonable situation all
$t$-structures associated with stability conditions are close to Noetherian (see Example in section
\ref{close-to-N-sec}).
However, it is important to observe that in the
non-Noetherian case the constant $t$-structures will usually lack some important features established
in \cite{AP} (such as the {\it open heart property}, see Proposition \ref{open-heart-prop}).

Next, we develop a little bit further the techniques of sheaves of $t$-structures by considering
an arbitrary morphism $f:X\to S$ and defining $t$-structures on $D(X)$, local over $S$. In the case
of a flat morphism we are able to define pull-backs of such $t$-structures under finite base changes of finite Tor dimension (see Theorem \ref{base-change-thm}). As a corollary we show that if $X$ is smooth then local (over $X$) $t$-structures on $D(X)$ are exactly the perverse $t$-structures constructed in
\cite{Bezr} (see Corollary \ref{local-cor}). Another application of this technique gives a description of the
heart of the constant $t$-structure on $D(X\times S)$ for affine scheme $S=\Spec(A)$ in terms of
$A$-modules in the heart of the corresponding $t$-structure on $D_{qc}(X)$ (see Proposition \ref{heart-prop}). 
We also show that if $L$ is an ample line bundle on $S$ then for a $t$-structure
on $X$ such that $f^*L\otimes D^{\le 0}(X)\sub D^{\le 0}(X)$ there exists
a new $t$-structure $(D_f^{\le 0}(X),D_f^{\ge 0}(X))$, local over $S$, with 
$D^{\ge 0}_f(X)=\cap_{n\ge 0} f^*L^{-n}\otimes D^{\ge 0}(X)$
(see Theorem \ref{localization-thm}).

Finally, we present one application of constant $t$-structures that seems to underscore once again the role of the Noetherian property. Namely, we prove that every bounded nondegenerate
Noetherian $t$-structure on $D(X)$ is invariant under the action of a connected group of autoequivalences of $D(X)$ (under a certain natural assumption on this action).

\noindent{\it Acknowledgment}. I am grateful to Dan Abramovich for his comments on the first
draft of the paper.

\noindent {\it Notation}. All our schemes are always assumed to be Noetherian of finite Krull dimension.
Starting from section \ref{AP-sec} they are
assumed to be of finite type over a fixed field $k$. We denote by $D(X)$ the bounded derived category of coherent sheaves on a scheme $X$, and by
$D_{qc}(X)$ the unbounded derived category of quasicoherent sheaves on $X$.
We denote the derived functor of tensoring simply by $\otimes$.
For a morphism of schemes $f:X\to Y$ we denote by $f_*$ and $f^*$ the derived functors of the
push-forward and the pull-back, respectively. In the case of a locally closed embedding $i:Y\hra X$ we
also use notation $F|_Y=i^*F$.
If $\AA\sub\CC$ is a subcategory in an additive category $\CC$ and $X\in\CC$ is an object then 
we write $\Hom(X,\AA)=0$ (resp., $\Hom(\AA,X)=0$) if $\Hom(X,Y)=0$ (resp., $\Hom(Y,X)=0$)
for all $Y\in\AA$. We define {\it left and right orthogonals} to $\AA$ in $\CC$ as the full subcategories given by
$\sideset{^{\perp}}{}{\AA}=\{X\in\CC \ | \Hom(X,\AA)=0\}$ and
$\AA^{\perp}=\{X\in\CC \ |\ \Hom(\AA,X)=0\}$, respectively.

\section{$t$-structures that are close to Noetherian ones}

\subsection{Preliminary remarks on $t$-structures and tiltings}

Our main reference for the theory of $t$-structures is section~1 of \cite{BBD}. Below we recall 
some basic definitions. 

Let $\DD$ be a triangulated category. 
A {\it $t$-structure} on $\DD$ is a pair of full subcategories
$(\DD^{\le 0},\DD^{\ge 0})$ satisfying the conditions (i) and (ii) below.
We denote $\DD^{\le n}=\DD^{\le 0}[-n]$, $\DD^{\ge n}=\DD^{\ge 0}[-n]$ for every $n\in\Z$.
Then the conditions are:

\noindent
(i) $\Hom(X,Y)=0$ for every $X\in \DD^{\le 0}$ and $Y\in\DD^{\ge 1}$;

\noindent
(ii) every object $X\in\DD$ fits into an exact triangle 
$$\tau^{\le 0}X\to X\to \tau^{\ge 1}X\to\ldots$$
with $\tau^{\le 0}X\in\DD^{\le 0}$, $\tau^{\ge 1}X\in\DD^{\ge 1}X$.

It is easy to see that $\DD^{\ge 1}$ is exactly the right orthogonal of $\DD^{\le 0}$
(resp., $\DD^{\le 0}$ is the left orthogonal of $\DD^{\ge 1}$), and
the terms of the above triangle are determined functorially (due to condition (i)).
Similarly, one defines other {\it truncation functors} $\tau^{\le n}$, $\tau^{\ge n}$ for $n\in\Z$.
The {\it heart} of the $t$-structure is $\CC=\DD^{\le 0}\cap\DD^{\ge 0}$. It is an abelian category.
The associated cohomology functors are defined by $H^0=\tau^{\le 0}\tau^{\ge 0}$,
$H^i(X)=H^0(X[i])$.
We will also use the notation
$\DD^{[a,b]}=\DD^{\ge a}\cap\DD^{\le b}$, where $[a,b]\sub\Z$ 
is a (possibly infinite on one side) interval.

Following \cite{AP} we will say that a $t$-structure is {\it nondegenerate} if $\cap_n\DD^{\le n}=\cap_n\DD^{\ge n}=0$ and $\cup_n\DD^{\le n}=\cup_n\DD^{\ge n}=\DD$. Note that this terminology is
not standard---in \cite{BBD} such a $t$-structure is called {\it bounded and nondegenerate}.

Let $\DD_1$ and $\DD_2$ be a pair of triangulated categories equipped with $t$-structures.
An exact functor $\Phi:\DD_1\to\DD_2$ is called {\it left (resp., right) $t$-exact} if
$\Phi(\DD_1^{\ge 0})\sub\DD_2^{\ge 0}$ (resp., $\Phi(\DD_1^{\le 0})\sub\DD_2^{\le 0}$).
A {\it $t$-exact functor} is a functor that is both left and right $t$-exact.
We will use later the following simple observation.

\begin{lem}\label{exact-fun-lem} Let $\DD_1$ and $\DD_2$ be a pair
of triangulated categories equipped with $t$-structures. 

\noindent
(i) Let $\Phi:\DD_1\to\DD_2$ be a $t$-exact functor with $\ker\Phi=0$, i.e.,
for any $F\in\DD_1$ such that $\Phi(F)=0$ one has $F=0$. Then for any interval
$[a,b]$ (possibly infinite on one side) one has
$$\DD_1^{[a,b]}=\{F\in\DD_1\ |\ \Phi(F)\in\DD_2^{[a,b]}\}.$$

\noindent
(ii) Let $(\Phi_n:\DD_1\to\DD_2)_{n\in\Z}$ be a family of exact functors such that 
for every $F\in\DD_1^{[a,b]}$ there exists an integer $N$ such that $\Phi_n(F)\in\DD_2^{[a,b]}$ 
for $n>N$.
Assume also that for any $F\in\DD_1$ such that $\Phi_n(F)=0$ for $n\gg 0$ one has $F=0$.
Then
$$\DD_1^{[a,b]}=\{F\in\DD_1\ |\ \Phi_n(F)\in\DD_2^{[a,b]}\text{ for }n\gg 0\}.$$
\end{lem}

\Pf . (i) Let us check that for $F\in\DD_1$ such that
$\Phi(F)\in\DD_2^{\le 0}$, one has $F\in\DD_1^{\le 0}$. Indeed, by $t$-exactness of $\Phi$
we have
$$\Phi(\tau^{\ge 1}F)=\tau^{\ge 1}\Phi(F)=0.$$
Hence, $\tau^{\ge 1}F=0$ by our assumption that $\ker\Phi=0$. 
It follows that $F\in\DD_1^{\le 0}$. Similarly, if $\Phi(F)\in\DD_2^{\ge 0}$ then $F\in\DD_1^{\ge 0}$.

\noindent
(ii) The proof is completely analogous to (i).
\ed

We refer to \cite{HRS} for basic facts about tilting with respect to a torsion theory.
Let $\CC$ be an abelian category. Recall that
a {\it torsion pair} $(\TT,\FF)$ in $\CC$ consists of two full subcategories
such that $\Hom(T,F)=0$ for every $T\in\TT$, $F\in\FF$, and such that every object $X\in\CC$ fits into
an exact sequence
$$0\to T\to X\to F\to 0$$
with $T\in\TT$, $F\in\FF$. 

Now assume that we have a $t$-structure $(\DD^{\le 0},\DD^{\ge 0})$ on a triangulated category $\DD$
and a torsion pair $(\TT,\FF)$ in the heart $\CC=\DD^{\le 0}\cap\DD^{\ge 0}$. Then one can define
a new $t$-structure $(\DD^{\le 0}_t,\DD^{\ge 0}_t)$ on $\DD$ by setting
$$\DD^{\le 0}_t=\{X\in\DD^{\le 1}\ |\ H^1X\in\TT\},
$$
$$\DD^{\ge 0}_t=\{X\in\DD^{\ge 0}\ |\ H^0X\in\FF\}.
$$
We say that this $t$-structure is obtained from $(\DD^{\le 0},\DD^{\ge 0})$ by {\it tilting with respect to
the torsion pair} $(\TT,\FF)$. The tilted heart $\CC_t=\DD^{\le 0}_t\cap\DD^{\ge 0}_t$
is equipped with a torsion pair $(\FF,\TT[-1])$.
Moreover, performing tilting with respect to this torsion pair will bring us back to the original
$t$-structure (up to a shift).\footnote{Since we do not require $\DD$ to be equivalent to the derived category of $\CC$, the new heart $\CC_t$ does not have to be equivalent to the abelian category
obtained by the tilting in the derived category of $\CC$, cf. Example 3.7 of \cite{Bridge2}.}
Note that we have $\DD^{\le 0}\sub\DD^{\le 0}_t\sub\DD^{\le 1}$. 
The following lemma shows
that this property characterizes pairs of $t$-structures related by tilting.

\begin{lem}\label{tilt-lem}
Let $(\DD_i^{\le 0},\DD_i^{\ge 0})$, $i=1,2$ be a pair of $t$-structures such that
\begin{equation}\label{tilt-in-eq}
\DD_1^{\le 0}\sub\DD_2^{\le 0}\sub\DD_1^{\le 1}
\end{equation}
Let us denote by $\CC_i$ the heart of $(\DD_i^{\le 0},\DD_i^{\ge 0})$ for $i=1,2$.
Then $(\DD_2^{\le 0},\DD_2^{\ge 0})$ is obtained from $(\DD_1^{\le 0},\DD_1^{\ge 0})$ by
tilting with respect to the torsion pair
$(\CC_2[1]\cap\CC_1,\CC_2\cap\CC_1)$ in $\CC_1$.
\end{lem}

\Pf . By passing to right orthogonals in \eqref{tilt-in-eq} (shifted by $[1]$) we find that
$$\DD_1^{\ge 1}\sub\DD_2^{\ge 0}\sub\DD_1^{\ge 0}.$$
Hence, $\CC_1\sub\DD_2^{[-1,0]}$ and $\CC_2\sub\DD_1^{[0,1]}$. 
Let us denote by
$\tau^*_i$ (resp., $H^*_i$) 
the truncation (resp., cohomology) functors associated with $(\DD_i^{\le 0},\DD_i^{\ge 0})$
for $i=1,2$. For any $X\in\CC_1$ consider the exact triangle
\begin{equation}\label{triangle}
A=\tau^{\le -1}_2X\to X\to \tau^{\ge 0}_2X=B\to\ldots
\end{equation}
Then $A\in\CC_2[1]$ and $B\in\CC_2$. Therefore, we have $H^i_1A=0$ for $i\neq -1,0$ and
$H^i_1B=0$ for $i\neq 0,1$. The long exact cohomology sequence associated with exact triangle \eqref{triangle} shows that $H^{-1}_1A=0$ and $H^1_1B=0$. 
Hence, both $A$ and $B$ belong to $\CC_1$. 
This proves that $(\CC_2[1]\cap\CC_1, \CC_2\cap\CC_1)$ is a torsion pair in $\CC_1$.
Switching the roles of $\CC_1[-1]$ and $\CC_2$ we derive that
any object $Y\in\CC_2$ fits into an exact triangle
$$A\to Y\to B\to A[1]$$
where $A=\tau_1^{\le 0}Y=H^0_1Y\in\CC_1\cap\CC_2$ and 
$B=\tau_2^{\ge 1}Y=(H^1_1Y)[-1]\in\CC_1[-1]\cap\CC_2$.
This implies that
$$\CC_2=\{Y\in\DD_1^{[0,1]}\ |\ H^0_1Y\in\CC_2\cap\CC_1, H^1_1Y\in\CC_2[1]\cap\CC_1\}$$
as required.
\ed

It is especially easy to construct torsion pairs in Noetherian abelian categories because of the
following simple observation.

\begin{lem}\label{Noeth-tor-lem} 
Let $\CC$ be a Noetherian abelian category. Then any full subcategory 
$\TT\sub\CC$ closed under quotients and extensions is contained in a torsion
pair $(\TT,\FF)$.
\end{lem}

\Pf . For every object $X\in\CC$ there is a unique maximal subobject of $X$ that belongs to $\TT$.
\ed

\begin{ex} If we have an increasing chain $\TT_1\sub\TT_2\sub\ldots $ of full subcategories closed under
quotients and extensions 
then the same is true for $\TT=\cup_n\TT_n$.
\end{ex}

\subsection{Pre-aisles that are close to Noetherian aisles}\label{close-to-N-sec}

Recall (see \cite{AJS})
that a full subcategory $\PP\sub\DD$ is called a {\it pre-aisle} if $\PP$ is closed under extensions and the shift functor $X\to X[1]$ (but not with respect to $X\to X[-1]$). 
A subcategory $\PP\sub\DD$ is called an {\it aisle} if $\PP=\DD^{\le 0}$ for some
$t$-structure on $\DD$. Clearly, every aisle is a pre-aisle. The converse is not true in general
(see Remark after Theorem \ref{induced-thm} below).

For a collection of subcategories $S_1,\ldots,S_n\sub\DD$
we denote by $\pra[S_1,\ldots,S_n]$ the 
smallest pre-aisle containing all $S_i$'s.
We call it the {\it pre-aisle generated by} $S_1,\ldots,S_n$.

\begin{defi} We say that a $t$-structure (or the corresponding aisle) is {\it Noetherian} if
its heart is Noetherian.
\end{defi} 
 
\begin{thm}\label{tilt-thm} 
Let $(\DD_0^{\le 0}, \DD_0^{\ge 0})$ be a Noetherian $t$-structure on $\DD$.
Then any pre-aisle $\PP\sub\DD$ such that $\DD^{\le -1}_0\sub\PP\sub\DD^{\le 0}_0$,
is an aisle, i.e., $\PP=\DD^{\le 0}$ for some $t$-structure on $\DD$.
\end{thm}

\Pf . Consider the heart $\CC_0=\DD_0^{\le 0}\cap\DD_0^{\ge 0}$.
Set $\TT=\CC_0\cap\PP$. Clearly, $\TT$ is stable under extensions. 
We claim that $\TT$ is also stable under taking quotients. Indeed, let 
$$X\to Y\to Z\to X[1]$$
be an exact triangle with $X,Y,Z\in\CC_0$ and with $Y\in\TT$. Then $Y\in\PP$ and
$$X[1]\in\DD_0^{\le 0}[1]=\DD_0^{\le -1}\sub\PP.$$
Hence, $Z\in\PP$ and therefore $Z\in\TT$.
Since $\CC_0$ is Noetherian, by Lemma \ref{Noeth-tor-lem} 
$\TT$ extends to a torsion pair $(\TT,\FF)$. We claim that 
$$\PP=\pra[\DD_0^{\le -1},\TT]=\{X\in\DD_0^{\le 0}\ |\ H^0X\in\TT\},$$
where $H^0$ is taken with respect to the $t$-structure $(\DD_0^{\le 0},\DD_0^{\ge 0})$.
Indeed, it suffices to check that for $X\in\PP$ one has $H^0X\in\TT$.
Consider the exact triangle
$$\tau^{\le -1}X\to X\to H^0X\to\tau^{\le -1}X[1].$$
We have $\tau^{\le -1}X[1]\in \DD_0^{\le -1}[1]=\DD_0^{\le -2}\sub\PP$. Therefore,
$H^0X\in\PP$ and hence $H^0X\in\TT$.
Thus, $\PP$ coincides with the aisle $\DD_t^{\le -1}$ 
of the tilted $t$-structure associated with $(\TT,\FF)$.
\ed

\begin{cor}\label{lim-cor} 
Let $(\DD_n^{\le 0},\DD_n^{\ge 0})_{n\ge0}$ be a sequence of $t$-structures in $\DD$ such that
$$\DD_0^{\le 0}\sub\DD_1^{\le 0}\sub\DD_2^{\le 0}\sub\ldots\sub\DD_0^{\le 1}.$$
Assume in addition that $\DD_0^{\le 0}\cap\DD_0^{\ge 0}$ is Noetherian. 
Then there exists a $t$-structure 
$(\DD_{\infty}^{\le 0},\DD_{\infty}^{\ge 0})$ with $\DD_{\infty}^{\le 0}=\cup_{n}\DD_n^{\le 0}$.
\end{cor}

\begin{defi} We say that a $t$-structure $(\DD^{\le 0},\DD^{\ge 0})$ on a triangulated category $\DD$
is {\it close to Noetherian} if there exists a Noetherian $t$-structure $(\DD_0^{\le 0},\DD_0^{\ge 0})$
on $\DD$
such that $\DD_0^{\le -1}\sub \DD^{\le 0}\sub \DD_0^{\le 0}$.
\end{defi}

In other words, close to Noetherian $t$-structures are precisely $t$-structures obtained
by tilting from Noetherian $t$-structures. 

\begin{ex} Let $\DD$ be a numerically finite triangulated category (see \cite{Bridge1}, sec. 1.3).
With every connected component $\Sigma$ of the space of numerical stability conditions
Bridgeland~\cite{Bridge1} associates a subspace $V(\Sigma)\sub(\NN(\DD)\otimes\C)^*$ such that
the map sending a stability to its central charge gives a local homeomorphism $\Sigma\to V(\Sigma)$.
Assume that $V(\Sigma)$ is defined over $\Q$ (this is true in all the known examples). 
Then for a dense subset $\Sigma_{\Q}\sub\Sigma$ the central charge $Z$ has the image in $\Q+i\Q$.
By Proposition 5.0.1 of \cite{AP} for $(\PP,Z)\in\Sigma_{\Q}$
the abelian category $\PP(t,t+1]$ will be Noetherian for a dense set of $t\in\R$.
It follows that for every stability $(\PP,Z)\in\Sigma$ the corresponding $t$-structure 
$(\PP(0,+\infty),\PP(-\infty,1])$ is close to Noetherian. Indeed, if $(\PP,Z)$ is sufficiently close
to $(\PP',Z')\in\Sigma_{\Q}$ then $\PP(0,1]\sub\PP'(-\eps,1+\eps]\sub\PP'(t,t+2]$ for
some $t\in\R$ such that $\PP'(t,t+1]$ is Noetherian.
\end{ex} 
 
\section{Induced $t$-structures}

\subsection{Abstract setting}\label{abs-sec}

Let $\TT$ be a triangulated category in which all small coproducts exist. 
Recall that a subcategory of $\TT$ is called
{\it cocomplete} if it is closed under small coproducts. 

\begin{defi}
For a subcategory $S\sub\TT$ we define the {\it cocomplete
pre-aisle generated by $S$}, denoted by $\pra[[S]]=\pra_{\TT}[[S]]$, as the smallest cocomplete
pre-aisle containing $S$. 
\end{defi}

Below we are going to use the following powerful theorem (Theorem A.1 of \cite{AJS}):
{\it for a small subcategory $S\sub\DD$ the pre-aisle $\pra[[S]]$ 
is an aisle}, i.e., $\pra[[S]]=\TT^{\le 0}$ for some $t$-structure $(\TT^{\le 0},\TT^{\ge 0})$
on $\TT$. Here is the first immediate application.

\begin{lem}\label{cocomplete-lem}
Let $\wt{\DD}$ be a triangulated category in which all small coproducts exist, and
let $\DD\sub\wt{\DD}$ be a full triangulated essentially small subcategory.
Then for every $t$-structure $(\DD^{\le 0},\DD^{\ge 0})$ on $\DD$ there exists a $t$-structure
$(\wt{\DD}^{\le 0},\wt{\DD}^{\ge 0})$ on $\wt{\DD}$ such that
$$\wt{\DD}^{\le 0}=\pra_{\wt{\DD}}[[\DD^{\le 0}]],$$
$$\wt{\DD}^{\ge 0}=\{F\in\wt{\DD}\ |\ \Hom(\DD^{\le -1},F)=0\}.$$
Furthermore, for every interval $[a,b]$ (possibly infinite on one side)
one has $\wt{\DD}^{[a,b]}\cap\DD=\DD^{[a,b]}$.
\end{lem}

\Pf . By Theorem A.1 of \cite{AJS} quoted above
there exists a $t$-structure $(\wt{\DD}^{\le 0},\wt{\DD}^{\ge 0})$
on $\wt{\DD}$ with $\wt{\DD}^{\le 0}=\pra_{\wt{\DD}}[[\DD^{\le 0}]]$. The formula for
$\wt{\DD}^{\ge 0}$ is easy to deduce (cf. Lemma 3.1 of \cite{AJS}). Note that since
$\DD$ is a full subcategory in $\wt{\DD}$, we have $\DD^{\ge 0}\sub\wt{\DD}^{\ge 0}$.
Therefore, the inclusion functor
$\DD\to\wt{\DD}$ is $t$-exact, and the last assertion
follows from Lemma \ref{exact-fun-lem}(i).
\ed

Let $\wt{\DD}_1$ and $\wt{\DD}_2$ be a pair of triangulated categories in which all small coproducts
exist, and let $\DD_1\sub\wt{\DD}_1$ and $\DD_2\sub\wt{\DD}_2$ be full triangulated essentially small subcategories.
Assume we have an exact functor $\Phi:\wt{\DD}_1\to\wt{\DD}_2$ commuting with small coproducts that admits a left adjoint functor $\Psi:\wt{\DD}_2\to\wt{\DD}_1$. Assume that $\Phi(\DD_1)\sub\DD_2$ and 
$\Psi(\DD_2)\sub\DD_1$.

\begin{thm}\label{induced-thm} 
(i) Let $(\DD_2^{\le 0},\DD_2^{\ge 0})$ be a $t$-structure on $\DD_2$ such that
the functor $\Phi\Psi:\DD_2\to\DD_2$ is right $t$-exact. Assume in addition that
$\DD_1=\Phi^{-1}(\DD_2)$, i.e., for any object $F\in\wt{\DD_1}$ such that $\Phi(F)\in\DD_2$ one has
$F\in\DD_1$. Then there exists a (unique) $t$-structure on $\DD_1$ with 
$$\DD_1^{\ge 0}=\{ F\in\DD_1\ |\ \Phi(F)\in\DD_2^{\ge 0}\}.$$
Moreover, the functor $\Phi$ is $t$-exact with respect to these $t$-structures.

\noindent
(ii) Assume also that for any $F\in\DD_1$ such that $\Phi(F)=0$ one has $F=0$. Then
$$\DD_1^{\le 0}=\{ F\in\DD_1\ |\ \Phi(F)\in\DD_2^{\le 0}\}.$$
In this situation if $\CC_2$ is Noetherian then so is $\CC_1$, where
$\CC_i=\DD_i^{\le 0}\cap\DD_i^{\ge 0}$.
\end{thm}

\Pf . (i) Let us extend the $t$-structure on $\DD_2$ to a $t$-structure 
$(\wt{\DD}_2^{\le 0},\wt{\DD}_2^{\ge 0})$ on $\wt{\DD}_2$ as in Lemma \ref{cocomplete-lem},
so that $\wt{\DD}_2^{\le 0}=\pra_{\wt{\DD}_2}[[\DD_2^{\le 0}]]$. 
Now let us define the $t$-structure on $\wt{\DD}_1$ by setting
$$\wt{\DD}_1^{\le 0}=\pra_{\wt{\DD}_1}[[\Psi(\DD_2^{\le 0})]]$$ 
(this is possible by Theorem A.1 of \cite{AJS}).
Then $\wt{\DD}_1^{\ge 1}$ is the right orthogonal to $\Psi(\DD_2^{\le 0})$ in $\wt{\DD}_1$.
Using adjointness of the pair $(\Psi,\Phi)$ we obtain that
$$\wt{\DD}_1^{\ge 0}=\{ F\in\wt{\DD}_1\ |\ \Phi(F)\in\wt{\DD}_2^{\ge 0}\}.$$
Now we claim that the functor $\Phi:\wt{\DD}_1\to\wt{\DD}_2$ is $t$-exact.
Indeed, clearly we have $\Phi(\wt{\DD}_1^{\ge 0})\sub\wt{\DD}_2^{\ge 0}$. Also, we have to check
that $\Phi(\wt{\DD}_1^{\le 0})\sub\wt{\DD}_2^{\le 0}$. Since $\Phi$ is exact and commutes with 
small coproducts, this follows from our assumption that $\Phi\Psi(\DD_2^{\le 0})\sub\DD_2^{\le 0}$.

Next we claim that setting $\DD_1^{[a,b]}=\wt{\DD}_1^{[a,b]}\cap\DD_1$ we get
a $t$-structure on $\DD_1$.
For this we need to prove that if $F\in\DD_1$ then $\tau^{\le 0}F\in\DD_1$ (apriori it lies in
$\wt{\DD}_1$). By $t$-exactness of $\Phi$ we get that 
$$\Phi\tau^{\le 0}F\simeq\tau^{\le 0}\Phi F\in\DD_2.$$
Since $\Phi^{-1}(\DD_2)=\DD_1$, this implies that $\tau^{\le 0}F\in\DD_1$.
Our formula for $\DD_1^{\ge 0}$ follows from the above formula for $\wt{\DD}_1^{\ge 0}$ and
from the fact that $\wt{\DD}_2^{\ge 0}\cap\DD_2=\DD_2^{\ge 0}$ (see Lemma \ref{cocomplete-lem}).

\noindent
(ii) The first assertion follows from Lemma \ref{exact-fun-lem}(i).
Now assume that $\CC_2$ is Noetherian.
Since $\Phi$ induces an exact functor with zero kernel
from $\CC_1=\DD_1^{\le 0}\cap\DD_1^{\ge 0}$ to $\CC_2$, this immediately implies
that $\CC_1$ is also Noetherian.
\ed


\begin{rem} According to Theorem A.1 of \cite{AJS} used above it is very easy to construct
$t$-structures in the unbounded derived category of quasicoherent sheaves $D_{qc}(X)$
by taking $D^{\le 0}_{qc}$ to be the cocomplete pre-aisle generated by some set of objects.
However, one should keep in mind that these $t$-structures rarely induce a $t$-structure on
$D(X)$. Here is the simplest example. Let $X$ be a smooth curve. Fix a point $p$
and define $D^{\le 0}_{qc}=\pra[[\OO_p]]$.
Then $D^{\ge 1}_{qc}$ consists of $F$ such that $\Hom^i(\OO_p,F)=0$ for $i\le 0$.
Consider the exact triangle in $D_{qc}(X)$
$$(j_*\OO_{X-p}/\OO_X)[-1]\to\OO_X\to j_*\OO_{X-p}\to j_*\OO_{X-p}/\OO_X,$$
where $j:X-p\to X$ is the natural open embedding.
It is easy to check that $j_*\OO_{X-p}\in D^{\ge 2}_{qc}$ while 
$(j_*\OO_{X-p}/\OO_X)[-1]\in D^{\le 1}_{qc}$. It follows that with respect to our $t$-structure one has
$\tau^{\ge 2}(\OO_X)=j_*\OO_{X-p}$, so the coherence is not preserved.
In other words, $D^{\le 0}_{qc}\cap D(X)$ is a pre-aisle, but not an aisle.
\end{rem}

\subsection{Applications to coherent sheaves: first examples}

The above abstract theorem can be applied in the case when $\DD_1$ and $\DD_2$ are bounded derived categories of coherent sheaves on some schemes and $\wt{\DD}_i$ are corresponding unbounded derived categories of quasicoherent sheaves (where $(\Psi,\Phi)$ is a
pair of adjoint exact functors of geometric origin). The simplest case when $\Phi$ is the push-forward
functor gives the following result.

\begin{prop}\label{finite-prop}
Let $f:X\to Y$ be a finite morphism of finite Tor dimension.
Assume that we have a $t$-structure $(D^{\le 0}(Y),D^{\ge 0}(Y))$ on $D(Y)$ such that
tensoring with $f_*\OO_X$ is a right $t$-exact functor.
Then there exists a $t$-structure on $D(X)$ with
$$D^{[a,b]}(X)=\{F\in D(X)\ |\ f_*F\in D^{[a,b]}(Y)\}$$
If $D^{\le 0}(Y)\cap D^{\ge 0}(Y)$ is Noetherian then so is
$D^{\le 0}(X)\cap D^{\ge 0}(X)$.
\end{prop}

\Pf . We simply have to apply Theorem \ref{induced-thm} to $\Phi=f_*$. Note that $f_*$ 
commutes with small coproducts (see Lemma 1.4 of \cite{Neeman}). The left adjoint is
$\Psi=f^*$, and we have $\Phi\Psi(F)\simeq f_*f^*F\simeq F\otimes f_*\OO_X$.
The assumption that $f$ has finite Tor dimension ensures that $f^*$ preserves boundedness of cohomology.
\ed

We have the following corollary for the theory of stability conditions (see \cite{Bridge1}).

\begin{cor}\label{stab-cor}
Let $f:X\to Y$ be a finite morphism of finite Tor dimension. Assume that we have a stability condition $(\PP,Z)$ on $D(Y)$ such that $f_*\OO_X\otimes\PP(t)\sub\PP(t,+\infty)$ for every $t\in\R$. Then
there exists an induced stability condition $(\PP',Z')$ on $D(X)$ with central charge $Z'=Z\circ f_*$ and
$\PP'(t)=\{F\ |\ f_*F\in\PP(t)\}$.
\end{cor}

\Pf . By Proposition \ref{finite-prop} we have a $t$-structure on $D(X)$ with the heart 
$\CC'=\PP'(0,1]=\{F\ |\ f_*F\in\PP(0,1]\}$ such that $Z'$ is the centered slope-function on $\CC'$.
Similarly, for every $t\in\R$ we can define a $t$-structure $(\PP'(>t),\PP'(\le t+1))$ on $D(X)$.
It remains to prove that the pair $(\CC',Z')$ satisfies the Harder-Narasimhan property.
Note that if $f_*F$ is semistable of phase $t\in(0,1]$ then so is $F$. 
Now for any $F\in\CC'$ let $G_1\sub G_2\sub\ldots\sub G_n=f_*F$ be the
Harder-Narasimhan filtration of $f_*F$, where $G_i/G_{i-1}$ is semistable of phase $t_i$.
Using the truncations with respect to the $t$-structures on $D(X)$ associated with $t_i$'s
we can construct a filtration $F_1\sub F_2\sub\ldots\sub F_n=f_*F$, such that $G_i=f_*F_i$.
Since $f_*(F_i/F_{i-1})=G_i/G_{i-1}$ is semistable of phase $t_i$, the same is true for
$F_i/F_{i-1}$.
\ed

\begin{ex} The assumptions of the above corollary are satisfied if
$f_*\OO_X=\oplus_i L_i$, where $L_i\in\Pic^0(Y)$, and the stability condition on $Y$
is stable under tensoring with $\Pic^0(Y)$. The latter condition can often be checked
(see Corollary \ref{invar-cor}).
\end{ex}

Here is another application of Theorem \ref{induced-thm}.

\begin{prop}\label{G-prop}
Let $G$ be a finite (discrete) group acting on $X$.
Then there is a bijection between $t$-structures on $D(X)$, 
invariant under $g^*:D(X)\to D(X)$ for every $g\in G$,
and $t$-structures on the derived category of equivariant coherent sheaves
$D_G(X)$ with respect to which the functor $D_G(X)\to D_G(X):F\mapsto F\otimes_G R$ is $t$-exact,
where $R$ is the regular representation of $G$.
Similarly, there is a bijection between stabilities on $D(X)$ invariant under $G$ and
stabilities $(\PP,Z)$ on $D_G(X)$ such that 
$\PP(t)\otimes_G R\sub\PP(t)$ for all $t\in\R$ and $Z(F\otimes_G R)=|G| Z(F)$ for all $F\in D_G(X)$.
\end{prop}

\Pf . Let $\Psi:D_G(X)\to D(X)$ denote the forgetful functor, and let
$\Phi:D(X)\to D_G(X)$ be the functor sending a coherent sheaf $F$ to 
the $G$-equivariant sheaf $\oplus_{g\in G} g^*F$.
Then both pairs $(\Phi,\Psi)$ and
$(\Psi,\Phi)$ are adjoint and we have natural isomorphisms
\begin{equation}\label{phipsi-eq}
\Phi\Psi F\simeq F\otimes_G R, \  \Psi\Phi F\simeq \oplus_{g\in G}g^*F,
\end{equation}
where $R$ is the regular representation of $G$. 
It remains to apply Theorem \ref{induced-thm} to both functors.

The bijection between stability conditions is established in a similar way.
If $(\PP,Z)$ is a $G$-invariant stability on $D(X)$ then we define a stability $(\PP',Z')$ on $D_G(X)$
by setting $\PP'(t)=\Psi^{-1}(\PP(t))$, $Z'=Z\circ\Psi$. We leave the details for the reader.
\ed


\subsection{Locality}\label{locality-sec}

The following general result appears as Corollary 2 in \cite{Bezr}. For completeness we supply
a proof (sketched in \cite{Bezr}).

\begin{lem}\label{open-ext} 
Let $X$ be a Noetherian scheme, $j:U\hra X$ an open subscheme. 
Then the restriction functor $j^*:D(X)\to D(U)$ is essentially surjective.
\end{lem}

\Pf . We use the following fact about extensions of sheaves (that follows easily from \cite{H}, Ex. II.5.15).
Let $F\to G$ be a surjective morphism of quasicoherent sheaves on a Noetherian scheme $X$ and
let $U\sub X$ be an open subset such that $G|_U$ is coherent. Then there exists a coherent
subsheaf $F'\sub F$ such that the induced morphism $F'|_U\to G|_U$ is surjective.

Now let $F^{\bullet}$ be a bounded complex of quasicoherent sheaves such that
$F^{\bullet}|_U$ has coherent cohomology. Let us denote $Z^i=\ker(d:F^i\to F^{i+1})$,
$B^i=d(F^{i-1})\sub Z^i$, $H^i=Z^i/B^i$.
Then using the above observation we can construct
a subcomplex of coherent sheaves $F_c^{\bullet}\sub F^{\bullet}$ such that 

\noindent
(i) $d(F_c^i)=(F^{i+1}_c\cap B^{i+1})$ and

\noindent
(ii) the natural map $(F^i_c\cap Z^i)|_U\to H^i|_U$ is surjective.

\noindent
It is easy to see that (i) and (ii) imply that $F_c^{\bullet}|_U$ is quasiisomorphic to $F^{\bullet}|_U$.
The subsheaves $F_c^i\sub F^i$ are constructed by descending induction in $i$.
Namely, assuming that $F_c^{i+1}$ is already constructed we first construct a coherent subsheaf
$Z_c^i\sub Z^i$ such that $Z_c^i|_U$ surjects onto $H^i|_U$ (by applying the above fact to the morphism $Z^i\to H^i$). This will guarantee condition (ii)
for any $F_c^i$ containing $Z_c^i$. Next, applying the above fact to the morphism
$d^{-1}(F^{i+1}_c)/Z_c^i\to F^{i+1}_c\cap B^{i+1}$ induced by $d$, we construct a subsheaf
$G_c\sub d^{-1}(F^{i+1}_c)/Z_c$ such that $d(G_c)=F^{i+1}_c\cap B^{i+1}$.
Finally, we let $F_c^i\sub d^{-1}(F^{i+1}_c)$ to be the preimage of $G_c$ under the natural
projection.

Now given $G\in D(U)$ we take $F=j_*G$, so that $F|_U\simeq G$, and construct $F_c\sub F$
as above. Then $F_c\in D(X)$ and $F_c|_U\simeq F|_U\simeq G$.
\ed



\begin{defi} Let $f:X\to S$ be a morphism of schemes. 
We say that a $t$-structure $(D^{\le 0}(X),D^{\ge 0}(X))$
on $D(X)$ is {\it local over } $S$
if for every open $U\sub S$ there exists a $t$-structure on $D(f^{-1}(U))$ such that the
restriction functor $D(X)\to D(f^{-1}(U))$ is $t$-exact.
A $t$-structure on $D(X)$ is called {\it local} if it is local over $X$ (with $f=\id$).
\end{defi}

By Lemma \ref{open-ext} the induced $t$-structure on $D(f^{-1}(U))$ is uniquely defined for
every open $U\sub S$.
It is easy to see that for $F\in D(X)$ the condition $F\in D^{[a,b]}(X)$ can be checked locally.
Indeed, since the cohomology functors $H^i$ with respect to our $t$-structures
commute with restrictions to open subsets, this follows immediately from the fact that the
condition $F=0$ for $F\in D(X)$ can be checked locally.
We derive also that for any vector bundle $V$ on $S$ the functor of tensoring with $f^*V$ is $t$-exact
with respect to a $t$-structure, local over $S$.
Finally, note that if a $t$-structure on $D(X)$ is nondegenerate and local over $S$ then the same is true for the induced $t$-structures on $D(f^{-1}(U))$ for any open $U\sub S$.

For example, $t$-structures on $D(X)$ considered in \cite{Bezr} (associated with monotone and comonotone perversity functions on the topological space of $X$) are local.
Below we will show that for smooth $X$ these are the only nondegenerate
local $t$-structures (see Corollary \ref{local-cor}).

The notion of a sheaf of $t$-structures considered in \cite{AP}, sec. 2.1, is equivalent to a $t$-structure
on $D(X\times S)$ local over $S$. Almost all asssertions made in {\it loc. cit.} about this notion
easily extend to the case of an arbitrary morphism $X\to S$.
The following theorem is a slight strengthening of Theorem 2.1.4 of \cite{AP}.

\begin{thm}\label{loc-thm} 
Let $f:X\to S$ be a morphism,
where $S$ is quasiprojective over an affine scheme. 
Let $L$ be an ample line bundle on $S$.
Then a nondegenerate $t$-structure $(D^{\le 0}(X),D^{\ge 0}(X))$ 
on $D(X)$ is local over $S$ iff tensoring with $f^*L$ is left $t$-exact, i.e.,
$f^*L\otimes D^{\ge 0}(X)\sub D^{\ge 0}(X)$.
\end{thm}

\Pf . The proof follows the same outline as that of Theorem 2.1.4 of \cite{AP}.
Let us observe that smoothness assumption used in {\it loc. cit.} is not necessary because of
Lemma \ref{open-ext}.
Another change is in the analogue of Lemma 2.1.6: we claim that it is enough
to assume only left $t$-exactness of tensoring with $f^*L$ to deduce that for a closed subset
$T\sub S$ an object $F\in D(X)$ is supported on $f^{-1}(T)$ iff all cohomology
objects $H^iF$ with respect to our $t$-structure are supported on $f^{-1}(T)$.
Indeed, one has to check that for an object $F\in D(X)$ and a section $s\in H^0(S,L^d)$ (where $d>0$) the vanishing of the morphism of multiplication by $s$
\begin{equation}\label{F-s-eq}
F\stackrel{s}{\ra} F\otimes f^*L^d
\end{equation}
induces the vanishing of the similar morphisms for
$\tau^{\ge 0}F$ and $\tau^{\le 0}F$. To this end consider the natural morphism
$$F\otimes f^*L^d\to\tau^{\ge 0}F\otimes f^*L^d.$$
By our assumption $\tau^{\ge 0}F\otimes f^*L^d\in D^{\ge 0}(X)$. Hence, the above morphism
factors through a morphism
$$\a:\tau^{\ge 0}(F\otimes f^*L^d)\to \tau^{\ge 0}F\otimes f^*L^d.$$
We claim that the composition of $\a$ with the morphism
$$\tau^{\ge 0}(F)\to\tau^{\ge 0}(F\otimes f^*L^d)$$
obtained  by applying $\tau^{\ge 0}$ to \eqref{F-s-eq}, coincides with the morphism
of multiplication by $s$ on $\tau^{\ge 0}F$. Indeed, it is enough to check this equality
after composing with the natural morphism $F\to\tau^{\ge 0}F$, so it follows from
the functoriality of the morphism \eqref{F-s-eq} in $F$.
\ed

\begin{cor} If $S$ is affine then any nondegenerate $t$-structure on $D(X)$ is local over $S$.
\end{cor}

There is a natural gluing procedure for $t$-structures, local over the base.

\begin{lem}\label{local-glue-lem} Let $f:X\to S$ be a morphism, and let $S=\cup_i U_i$ be
a finite open covering of $S$. Assume that for every $i$ we have a nondegenerate $t$-structure
on $D(f^{-1}(U_i))$, local over $U_i$, and that these $t$-structures agree over all pairwise
intersection. Then there exists a $t$-structure on $D(X)$, local over $S$, inducing
the given $t$-structure on every $D(f^{-1}(U_i))$. 
\end{lem}

\Pf . Let us set $X_i=f^{-1}(U_i)$. We want to check that
$$D^{[a,b]}(X)=\{ F \ |\ F|_{X_i}\in D^{[a,b]}(X_i)\text{ for all }i\}$$ 
is a $t$-structure on $X$. To show orthogonality of 
$F\in D^{\le 0}(X)$ and $G\in D^{\ge 1}(X)$
consider the object $R\HHom_S(F,G):=f_*R\HHom(F,G)\in D_{qc}(S)$.
Note that for every open subset $U\sub S$ we
have $R\Hom(F|_{f^{-1}U},G|_{f^{-1}U})\simeq R\Ga(U,R\HHom_S(F,G))$. Hence, for
every $i$ and every open affine $U\sub U_i$ we have 
$R\HHom_S(F,G)|_U\in D^{\ge 1}(U)$ (with respect to the standard
$t$-structure on $D(U)$). It follows that $R\HHom_S(F,G)\in D^{\ge 1}(S)$ and therefore
$\Hom(F,G)=0$. To define the truncation functors it suffices to define
$H^0F$ and $\tau^{\ge 1}F$ for $F\in D^{\ge 0}(X)$ (since the $t$-structures
on $D(X_i)$ are nondegenerate). Set $F_i=F|_{X_i}$. Then we have a natural gluing
datum for the objects $H^0F_i$ in the hearts of $t$-structures on $X_i$. 
At this point we observe that analogues of Theorem 2.1.9 and of
Corollary 2.1.11 of \cite{AP} hold in the situation of a general morphism $X\to S$ (with the same
proof). Therefore, we can glue $(H^0F_i)$ into an object $H^0F$ equipped with isomorphisms
$H^0F|_{X_i}\simeq H^0F_i$. Looking at restrictions to $X_i$'s one easily checks
that $R\HHom_S(H^0F,F)\in D^{\ge 0}(S)$. By the analogue of Lemma 2.1.10 of \cite{AP} this implies that we can glue morphisms $H^0F_i\to F_i$ into a global morphism $H^0F\to F$.
\ed

Using Proposition \ref{finite-prop} we get the following base change construction.

\begin{thm}\label{base-change-thm}
Let $f:X\to S$ be a flat morphism, and let $(D^{\le 0}(X), D^{\ge 0}(X))$ be a $t$-structure on $D(X)$, local over $S$. Then for any finite morphism of finite Tor dimension $g:S'\to S$ there is an induced
$t$-structure on $D(X\times_S S')$ given by
$$D^{[a,b]}(X\times_S S')=\{F\in D(X\times_S S')\ |\ g'_*(F)\in D^{[a,b]}(X)\},$$
where $g':X\times_S S'\to X$ is the natural projection. If the original $t$-structure on $D(X)$ is
Noetherian then so is the induced $t$-structure on $D(X\times_S S')$.
\end{thm}

\Pf .  We claim that in this case the assumptions of Proposition \ref{finite-prop}
are satisfied for the morphism $g':X\times_S S'\to X$. Indeed, it suffices to check
that tensoring with $g'_*\OO_{X\times_S S'}$ is right $t$-exact.
The question is local over $S$, so we can assume that $g_*\OO_{S'}$ has
a finite resolution $V_n\to\ldots\to V_0$ by vector bundles on $S$ (since $g$ has a finite Tor dimension). 
Then in the derived category
$D(S)$ we have
$$g_*\OO_{S'}\simeq (V_n\to\ldots\to V_0),$$
where the complex is concentrated in degrees $[-n,0]$.
Using the base change formula we get
$$g'_*\OO_{X\times_S S'}\simeq f^*g_*\OO_{S'}\simeq (f^*V_n\to\ldots\to f^*V_0).$$
Since our $t$-structure on $D(X)$ is local over $S$, the functors of tensoring
with $f^*V_i$ are $t$-exact. This implies that tensoring with
the above complex is right $t$-exact. It remains to apply Proposition \ref{finite-prop}.
\ed

For example, assume that $T\sub S$ is a closed subscheme that is a locally complete intersection. 
Then by the above theorem, a $t$-structure on $D(X)$, local over $S$, induces a $t$-structure on 
$D(f^{-1}(T))$, local over $T$, such that the push-forward functor $D(f^{-1}(T))\to D(X)$ is $t$-exact.

\begin{cor}\label{local-cor}
Assume that $X$ is smooth over a field $k$. Then any nondegenerate
local $t$-structure on $D(X)$ is the $t$-structure
associated with a monotone and comonotone perversity function on the topological
space of $X$ (see \cite{Bezr}).
\end{cor}

\Pf . By Theorem \ref{base-change-thm} (applied to $f=\id:X\to X$)
for every closed subscheme $i:Z\hra X$ we have an induced local 
$t$-structure on $D(Z)$ such that the functor $i_*:D(X)\to D(X)$ is $t$-exact. 
Furthermore, it is easy to see that this $t$-structure on $D(Z)$ is nondegenerate.
Now assume that $Z$ is irreducible and
reduced. By locality the $t$-structure on $D(Z)$ induces a nondegenerate $t$-structure on $D(\eta_Z)$,
where $\eta_Z\in Z$ is the generic point, such that the restriction functor $D(Z)\to D(\eta_Z)$
is $t$-exact. There is a unique integer $p=p(\eta_Z)$ such that this $t$-structure on $D(\eta_Z)$
coincides with $(D^{\le p}_{st}(\eta_Z),D^{\ge p}_{st}(\eta_Z))$, where $(D^{\le 0}_{st},D^{\ge 0}_{st})$
denotes the standard $t$-structure. Indeed, since $D(\eta_Z)$ is semisimple (i.e., every exact triangle splits) and the object
$k(\eta_Z)\in D(\eta_Z)$ is indecomposable, it has a unique nonzero cohomology with respect
to any nondegenerate $t$-structure.

Next, let us show that if we have an embedding
of irreducible closed subsets $Z\sub Y$ then
$p(\eta_Z)\ge p(\eta_Y)$, i.e., the function $p$ is monotone. 
By locality it suffices to study the situation in a neighborhood of
$\eta_Z$ in $Y$. Let $A$ be the local ring of $Y$ at $\eta_Z$, $\CC$ the heart of the induced
$t$-structure on $D(\Spec A)$, and let $k(\eta_Y)$, $k(\eta_Z)$ denote the residue fields at $\eta_Y$ and $\eta_Z$. Since $k(\eta_Y)[-p(\eta_Y)]$ belongs to the heart of the $t$-structure on $D(\eta_Y)$, there exists $F\in\CC$ such that $F|_{\eta_Y}\simeq k(\eta_Y)[-p(\eta_Y)]$.
On the other hand, viewing $k(\eta_Z)$ as an $A$-module we have $k(\eta_Z)[-p(\eta_Z)]\in\CC$.
Assume that $M=H_{st}^nF\neq 0$ and $H_{st}^{>n}F=0$, where $H_{st}^i$ denote the cohomology
functors with respect to the standard $t$-structure on $D(\Spec A)$. Note that $n\ge p(\eta_Y)$.
By Nakayama lemma $M$ has a nonzero morphism to $k(\eta_Z)$. Thus, we get a nonzero morphism
$F\to k(\eta_Z)[-n]$. Since $k(\eta_Z)[-n]\in\CC[p(\eta_Z)-n]$, this implies that $p(\eta_Z)-n\ge 0$.
Hence, $p(\eta_Z)\ge n\ge p(\eta_Y)$.

Let $\Du:D(X)\to D(X): F\mapsto R\und{\Hom}(F,\om_X[\dim X])$ be the duality functor. Then
$(\Du(D^{\ge 0}(X)),\Du(D^{\le 0})))$ is also a nondegenerate local $t$-structure on $D(X)$,
so we can apply the above argument to it as well. We claim that the  
corresponding function on points is $\ov{p}(x)=-\dim x-p(x)$. 
Indeed, for every irreducible closed subset $Z\sub X$ there exists an object $F$ in the heart of the 
original $t$-structure on $D(Z)$ such that $F|_{\eta_Z}\simeq k(\eta_Z)[-p(\eta_Z)]$. 
Then $\Du(F)|_{\eta_Z}\simeq k(\eta_Z)[\dim Z+p(\eta_Z)]$. But $\Du(F)$ is in the heart of
the new $t$-structure, so we obtain the above formula for $\ov{p}$.
Since $\ov{p}$ is monotone, we get
that $p$ is comonotone. Finally, using the fact that for every closed subset $i:Z\hra X$ the functor
$i^*$ (resp., $i^!$) is right $t$-exact (resp., left $t$-exact)
we easily see that $D^{\le 0}\sub D^{p,\le 0}$ (resp., $D^{\ge 0}\sub D^{p,\ge 0}$), where
$(D^{p,\le 0}, D^{p,\ge 0})$ is the $t$-structure associated with $p$ (see \cite{Bezr}).
Therefore, we have equalities $D^{\le 0}=D^{p,\le 0}$, $D^{\ge 0}=D^{p,\ge 0}$.
\ed

Exactly the same argument as in Proposition 3.3.2 of \cite{AP} proves the following
{\it open heart property} of a Noetherian $t$-structure on $D(X)$, local over $S$.

\begin{prop}\label{open-heart-prop}
Let $f:X\to S$ be a flat morphism, 
$(D^{\le 0}(X), D^{\ge 0}(X))$ a {\em Noetherian} $t$-structure on $D(X)$, local over $S$.
For every open subset $U\sub S$ we denote by $\CC_U\sub D(f^{-1}(U))$ the heart of the
corresponding $t$-structure.
Let also $T\sub S$ be a closed subscheme that is a locally complete intersection, and
let $\CC_T\sub D(f^{-1}(T))$ be the heart of the induced $t$-structure. 
Then for every $F\in D(X)$ such that $F|_{f^{-1}(T)}\in\CC_T$ there exists an open
neighborhood $T\sub U\sub S$ such that $F|_{f^{-1}(U)}\in\CC_U$.
\end{prop}

The Noetherian hypothesis is used in the proof to guarantee the existence for any object in $\CC_S$
of a maximal subobject supported on $f^{-1}(T)$. 
Without this hypothesis the result is false
(see the counterexample in \cite{AP}, sec.~6.3).
In section \ref{invar-sec} we will describe one nice application of the open heart property to the invariance of a $t$-structure with respect to a continuous group of autoequivalences. 

For a future reference we record the following technical result about the base change with respect
to a closed embedding.

\begin{lem}\label{closed-emb-lem}
Let $f:X\to S$ be a flat morphism, and let $(D^{\le 0}(X), D^{\ge 0}(X))$ be a $t$-structure on $D(X)$, local over $S$. Then for any closed embedding 
of finite Tor dimension $i:T\to S$ we have
$$H^0k^*k_*F\simeq F,$$
where $k:X_T:=f^{-1}(T)\to X$ is the natural embedding, $H^0$ is taken with respect to
the induced $t$-structure on $D(X_T)$ and $F$ is an object in the heart of this $t$-structure.
\end{lem}

\Pf . Let us define $G\in D(X_T)$ from the exact triangle
$$G\to k^*k_*F\stackrel{\a}{\ra} F\to\ldots,$$
where $\a$ is the natural adjunction morphism. 
It suffices to show that $G\in D^{\le -1}(X_T)$. By the definition of the $t$-structure on $D(X_T)$
this is equivalent to $k_*G\in D^{\le -1}(X)$. But in the exact triangle
$$k_*G\to k_*k^*k_*F\stackrel{k_*\a}{\ra} k_*F\to\ldots$$
the morphism $k_*\a$ is the projection onto the direct summand. Indeed, if
$\b:k_*F\to k_*k^*(k_*F)$ is the natural adjunction morphism then $k_*\a\circ\b=\id_{k_*F}$.
Hence, we have $k_*k^*k_*F\simeq k_*F\oplus k_*G$. On the other hand,
by the projection formula $k_*k^*k_*F\simeq k_*\OO_{X_T}\otimes k_*F$.
Moreover, the morphism $\b$ is induced by the natural map $\OO_X\to k_*\OO_{X_T}$
that has the cone $f^*J_T[1]$, where $J_T\sub\OO_S$ is the ideal sheaf of $T$.
Therefore,
$$k_*G\simeq f^*J_T\otimes k_*F[1].$$
Using locality of the $t$-structure and local finite resolutions of $J_T$ over $\OO_S$
we derive that $k_*G\in D^{\le -1}(X)$.
\ed

\section{Constant families of $t$-structures}

\subsection{Gluing of $t$-structures}
 
Let us recall some definitions and constructions involving admissible subcategories and
semiorthogonal decompositions (see \cite{BK}, \cite{BO}, \cite{Orlov2}). 

\begin{defi}
Let $\DD$ be a triangulated category.
A full triangulated subcategory $\AA\sub\DD$ is called {\it right admissible} (resp.,
{\it left admissible})
if there exist right (resp., left) adjoint functors to the inclusion $\AA\to\DD$. 
\end{defi}

\begin{defi}
A {\it (weak) semiorthogonal decomposition}
\begin{equation}\label{semi-eq}
\DD=\lan\AA_1,\ldots,\AA_n\ran
\end{equation}
is given by a collection of full triangulated subcategories 
such that $\Hom(\AA_j,\AA_i)=0$ for $i<j$ and for every object $X\in\DD$ there exists a sequence
of exact triangles 
$$A_i\to X_i\to X_{i-1}\to A_i[1], \ i=1,\ldots, n,$$
with $A_i\in\AA_i$, where $X_n=X$ and $X_0=0$.
\end{defi}

We will not use the stronger notion of {\it semiorthogonal decomposition} that requires
all subcategories $\AA_i$ to be right and left admissible, so we will omit the attribute "weak".

In the case $n=2$ the semiorthogonal decomposition is determined by one of the subcategories
$\AA_1$ and $\AA_2$. Namely, $\DD=\lan\AA_1,\AA_2\ran$ iff $\AA_1$ is left admissible
and $\AA_2=\sideset{^{\perp}}{}{\AA_1}$. Equivalently, $\AA_2$ should be right admissible and
$\AA_1=\AA_2^{\perp}$. In general if \eqref{semi-eq} is a semiorthogonal decomposition
then $\AA_1$ is left admissible and there is a semiorthogonal decomposition
\begin{equation}\label{semi-perp-eq}
\sideset{^{\perp}}{}{\AA_1}=\lan\AA_2,\ldots,\AA_n\ran.
\end{equation}
Also, $\AA_n$ is right admissible and
$$\AA_n^{\perp}=\lan\AA_1,\ldots,\AA_{n-1}\ran.$$
This leads to an alternative definition of a semiorthogonal decomposition 
as a sequence of left admissible subcategories $\DD_1=\AA_1\sub\DD_2\sub\ldots\sub\DD_n=\DD$
such that $\AA_i$ is the left orthogonal of $\DD_{i-1}$ in $\DD_i$.

In the next lemma we define the gluing of $t$-structures in the situation when one has a semiorthogonal decomposition. This is a particular case of the formalism
of \cite{BBD}, sec.~1.4, rewritten in slightly different terms. It will be convenient to use the following notion analogous to that of a pre-aisle.
We say that a full subcategory $\QQ\sub\DD$ is an {\it anti-pre-aisle} if it is closed under extensions
and under the functor $X\to X[-1]$. For a collection of subcategories $S_1,\ldots,S_n\sub\DD$ we
denote by $\apra[S_1,\ldots,S_n]$ the anti-pre-aisle generated by $S_1,\ldots,S_n$
(i.e., the smallest anti-pre-aisle containing all $S_i$'s).

\begin{lem}\label{glue-lem} 
Assume we are given a semiorthogonal decomposition \eqref{semi-eq}
and $t$-structures $(\AA^{\le 0}_i,\AA^{\ge 0}_i)$ on each $\AA_i$.
If all the subcategories $\AA_i\sub\DD$ are right admissible
then the following formulas define a $t$-structure on $\DD$:
 \begin{equation}
\begin{array}{l}
\DD^{\le 0}_{\rho}=\pra[\AA^{\le 0}_1,\ldots,\AA^{\le 0}_n],\\
\DD^{\ge 0}_{\rho}=\{X\in\DD \ |\ \rho_1(X)\in\AA^{\ge 0}_1,\ldots,\rho_n(X)\in\AA^{\ge 0}_n\},\\
\end{array}
\end{equation}
where for each $i$ the functor
$\rho_i:\DD\to\AA_i$ is the right adjoint to the inclusion $\AA_i\to\DD$.
Similarly, if all the subcategories $\AA_i\sub\DD$ are left admissible then
we have a $t$-structure on $\DD$ defined by
\begin{equation}
\begin{array}{l}
\DD^{\le 0}_{\la}=\{X\in\DD \ |\ \la_1(X)\in\AA^{\le 0}_1,\ldots,\la_n(X)\in\AA^{\le 0}_n\},\\
\DD^{\ge 0}_{\la}=\apra[\AA^{\ge 0}_1,\ldots,\AA^{\ge 0}_n],
\end{array}
\end{equation}
where $\la_i:\DD\to\AA_i$ are the left adjoint functors to the inclusions.
\end{lem}

\Pf . Let us prove that $(\DD_{\rho}^{\le 0}, \DD_{\rho}^{\ge 0})$ is a $t$-structure on $\DD$
provided all $\AA_i$'s are right admissible
(for the second $t$-structure the argument is similar). First, let us consider the case $n=2$.
To show orthogonality of
$\DD_{\rho}^{\le 0}$ and $\DD_{\rho}^{\ge 1}$ it is enough to check that
for $X\in\DD_{\rho}^{\ge 1}$ one has
$\Hom(\AA_i^{\le 0}, X)=0$ for $i=1,2$. 
But this follows immediately from the assumption $\rho_i X\in\AA_i^{\ge 1}$ and from the
orthogonality of $\AA_i^{\le 0}$ and $\AA_i^{\ge 1}$.
Next, for every $X\in\DD$ we have to find an exact triangle
$Y\to X\to Z\to Y[1]$ with $Y\in\DD_{\rho}^{\le 0}$ and $Z\in\DD_{\rho}^{\ge 1}$.
By the definition of $\rho_2$ we have an exact triangle
\begin{equation}\label{glue-tr1}
\rho_2(X)\to X\to W\to\rho_2(X)[1],
\end{equation}
where $\rho_2(W)=0$.
Also, we have an exact triangle
$$\tau^{\le 0}\rho_2(X)\to\rho_2(X)\to\tau^{\ge 1}\rho_2(X)\to\tau^{\le 0}\rho_2(X)[1],$$
where we use the truncation functors on $\AA_2$.
We can embed the composed map $\tau^{\le 0}\rho_2(X)\to\rho_2(X)\to X$ into an exact
triangle
$$\tau^{\le 0}\rho_2(X)\to X\to X'\to\tau^{\le 0}\rho_2(X)[1]$$
By the octahedron axiom we also have an exact triangle
$$\tau^{\ge 1}\rho_2(X)\to X'\to W\to\tau^{\ge 1}\rho_2(X)[1]$$
that implies that $\rho_2(X')\simeq\tau^{\ge 1}\rho_2(X)$. 
Similarly, we can embed the composed map
$\tau^{\le 0}\rho_1(X')\to\rho_1(X')\to X'$ (the truncation is taken in $\AA_1$)
into an exact triangle
$$\tau^{\le 0}\rho_1(X')\to X'\to Z\to\tau^{\le 0}\rho_1(X')[1],$$
where $\rho_1(Z)\simeq\tau^{\ge 1}\rho_1(X')$.
Also, since $\rho_2(\AA_1)=0$, it follows that
$\rho_2(Z)\simeq\rho_2(X')\simeq\tau^{\ge 1}\rho_2(X)$.
Therefore, $Z\in\DD_{\rho}^{\ge 1}$. Finally, let us embed
the composed map $X\to X'\to Z$ into an exact triangle
$Y\to X\to Z\to Y[1]$.
By the octahedron axiom we have an exact triangle
$$\tau^{\le 0}\rho_2(X)\to Y\to\tau^{\le 0}\rho_1(X')\to\tau^{\le 0}\rho_2(X)[1],$$
hence, $Y\in\DD_{\rho}^{\le 0}$.

The case of general $n$ is deduced by induction:
one has to use the semiorthogonal decompositions
\eqref{semi-perp-eq} and
$\DD=\lan\AA_1,\sideset{^{\perp}}{}{\AA_1}\ran$
(note that $\sideset{^{\perp}}{}{\AA_1}$ is automatically right admissible).
\ed

Under additional assumptions one can rewrite the above definition of the glued $t$-structure
in a more symmetric way. We keep the notation $\rho_i$ (resp., $\la_i$) for the right
(resp., left) adjoint functor to the inclusion $\AA_i\to\DD$.

\begin{lem}\label{glue-lem2} 
Assume we are given a semiorthogonal decomposition \eqref{semi-eq}
and $t$-structures $(\AA^{\le 0}_i,\AA^{\ge 0}_i)$ on each $\AA_i$.
Assume in addition that all the subcategories $\AA_i$ are right admissible,
and for every $i<j$ the functor $\rho_i|_{\AA_j}:\AA_j\to\AA_i$ is
right $t$-exact (with respect to the $t$-structures on $\AA_i$ and $\AA_j$). Then one has
\begin{equation}\label{simple-t-str-right}
\DD^{[a,b]}_{\rho}=\{X\in\DD\ |\ \rho_1(X)\in\AA_1^{[a,b]},\ldots,\rho_n(X)\in\AA_n^{[a,b]}\}.
\end{equation}
Similarly, if we assume that $\AA_i$'s are left admissible, and
for every $i<j$ the functor $\la_j|_{\AA_i}:\AA_i\to\AA_j$ is
left $t$-exact, then 
\begin{equation}\label{simple-t-str-left}
\DD^{[a,b]}_{\la}=\{X\in\DD\ |\ \la_1(X)\in\AA_1^{[a,b]},\ldots,\la_n(X)\in\AA_n^{[a,b]}\}.
\end{equation}
\end{lem}

\Pf . First, consider the case $n=2$. Note that 
for every $X\in\DD$ the exact triangle \eqref{glue-tr1} can be rewritten as
\begin{equation}\label{glue-tr2}
\rho_2(X)\to X\to\la_1(X)\to\rho_2(X)[1].
\end{equation}
Also, we have
$$\DD^{\le 0}_{\rho}=\{X\in\DD\ |\ \la_1(X)\in\AA_1^{\le 0}, \rho_2(X)\in\AA_2^{\le 0}\}.$$
Applying $\rho_1$ to \eqref{glue-tr2} we get the exact triangle
$$\rho_1\rho_2(X)\to\rho_1(X)\to\la_1(X)\to\rho_1\rho_2(X)[1].$$
Thus, if $\rho_2(X)\in\AA_2^{\le 0}$ then by our assumption $\rho_1\rho_2(X)\in\AA_1^{\le 0}$.
Thus, if $\rho_2(X)\in\AA_2^{\le 0}$ then the conditions $\rho_1(X)\in\AA_1^{\le 0}$ and
$\la_1(X)\in\AA_1^{\le 0}$ are equivalent. This proves that
$$\DD^{\le 0}_{\rho}=\{X\in\DD\ |\ \rho_1(X)\in\AA_1^{\le 0}, \rho_2(X)\in\AA_2^{\le 0}\}.$$

The case $n>2$ follows by induction. 
Namely, we apply the above argument to the semiorthogonal decompositions
\eqref{semi-perp-eq} and
$\DD=\lan\AA_1,\sideset{^{\perp}}{}{\AA_1}\ran$.
We only have to check that the restriction of the functor
$\rho_1$ to $\sideset{^{\perp}}{}{\AA_1}=\lan\AA_2,\ldots,\AA_n\ran$ is right $t$-exact
(with respect to the glued $t$-structure on $\lan\AA_2,\ldots,\AA_n\ran$).
But this immediately follows from the fact that
$\lan\AA_2,\ldots,\AA_n\ran^{\le 0}=\pra[\AA_2^{\le 0},\ldots,\AA_n^{\le 0}]$.
\ed


\subsection{Constant $t$-structure over $\P^r$ following \cite{AP}}\label{AP-sec}

Starting from this point 
we always assume our schemes to be of finite type over a fixed field $k$. The product of such schemes
is taken over $k$.

Let $X$ be a scheme. 

\begin{thm} (Theorem 2.3.6) of \cite{AP}). For every Noetherian nondegenerate
$t$-structure $(D^{\le 0}(X),D^{\ge 0}(X))$ on $D(X)$
there is a Noetherian nondegenerate $t$-structure
$(D^{\le 0}(X\times\P^r),D^{\ge 0}(X\times\P^r))$, local over $\P^r$, characterized by the property
\begin{equation}\label{const-t-str-eq}
D^{[a,b]}(X\times\P^r)=\{F\in D(X\times\P^r) \ |\ p_*(F(n))\in D^{[a,b]}(X)\text{ for all }n\gg 0\},
\end{equation}
where $p:X\times\P^r\to X$ is the natural projection. 
\end{thm}

We will refer to the above $t$-structure on $D(X\times\P^r)$ as the {\it constant $t$-structure}. 
We show in \cite{AP} that it is obtained as a certain limit of the sequence of glued $t$-structures on
$D(X\times\P^r)$. More precisely, one has
\begin{equation}\label{const-lim-eq}
D^{\le 0}(X\times\P^r)=\cup_n D^{\le 0}_n(X\times\P^r),
\end{equation}
where the $n$-th $t$-structure 
$(D_n^{\le 0}, D_n^{\ge 0})$ is glued from standard $t$-structures with respect to the
semiorthogonal decomposition
$$D(X\times\P^r)=\lan p^*D(X)(-r-n),\ldots,p^*D(X)(-1-n), p^*D(X)(-n)\ran.$$
In the notation of Lemma \ref{glue-lem} we have $D_n^{[a,b]}=D_{\rho}^{[a,b]}$.
The conditions of Lemma \ref{glue-lem2} are also satisfied in this case, so the $n$-th
$t$-structure has a nice description in terms of the functors $F\mapsto p_*(F(m))$ which leads to
\eqref{const-t-str-eq}. 

It is important to observe that the assumption that $X$ is smooth made in \cite{AP} is not needed for Theorem 2.3.6 (nor for Theorem 2.1.4) of {\it loc. cit.}. Indeed, it is used there only to guarantee the essential surjectivity of the restriction functor $j^*:D(X\times S)\to D(X\times U)$ for open subsets $U\sub S$. However, this is true without this assumption (see Lemma \ref{open-ext}).
The smoothness is used seriously in section 2.4 of {\it loc. cit.} to characterize the essential
image of the push-forward under a closed embedding. 
This characterization is then used in {\it loc. cit.} to construct the constant $t$-structure on $D(X\times S)$ for arbitrary smooth quasiprojective base $S$.
Using Theorem \ref{induced-thm} we will give an alternative construction of such a $t$-structure on
$D(X\times S)$ assuming only that $X$ and $S$ are of finite type of $k$.
We will also extend the construction of constant $t$-structures 
to the case of close to Noetherian $t$-structures (see section \ref{close-to-N-sec}).

\subsection{Constant $t$-structures}\label{const-t-str-sec}

We start with the case when the base is $\P^r$. 
Using Lemma \ref{glue-lem} we can make the following observation.

\begin{lem}\label{gener-lem} 
Let $(D^{\le 0}(X),D^{\ge 0}(X))$ be a
Noetherian nondegenerate $t$-structure  on $D(X)$.
Then the corresponding constant $t$-structure on $D(X\times\P^r)$ satisfies
$$D^{\le 0}(X\times\P^r)=\pra[D^{\le 0}(X)\boxtimes\OO_{\P^r}(n)\ |\ n\in\Z].$$
\end{lem}

\Pf . This immediately follows from \eqref{const-lim-eq} and the formula for $D_n^{\le 0}(X\times\P^r)$  obtained by Lemma \ref{glue-lem}.
\ed

Now we can consider the case of a close to Noetherian $t$-structure on $D(X)$.

\begin{lem}\label{proj-lem} 
Let $(D^{\le 0}(X),D^{\ge 0}(X))$ be a close to Noetherian nondegenerate $t$-structure
on $D(X)$. Then we can define a $t$-structure on $D(X\times\P^r)$, local over $\P^r$, by
the formula
\begin{equation}\label{proj-eq2}
D^{[a,b]}(X\times\P^r)
=\{F\in D(X\times \P^r)\ |\ p_*(F\otimes \OO(n))\in D^{[a,b]}(X)\text{ for all }n\gg 0\}.
\end{equation}
We also have
\begin{equation}\label{const-gener-eq}
D^{\le 0}(X\times\P^r)=\pra[D^{\le 0}(X)\boxtimes\OO(n)\ |\ n\in\Z].
\end{equation}
\end{lem}

\Pf . Suppose that $(D^{\le 0}_0(X), D^{\ge 0}(X))$ is a Noetherian $t$-structure such that
$$D^{\le -1}_0(X)\sub D^{\le 0}(X)\sub D^{\le 0}_0(X).$$ 
Then we have
the corresponding Noetherian constant
$t$-structure $(D^{\le 0}_0(X\times \P^r), D^{\ge 0}_0(X\times \P^r))$ on $D(X\times\P^r)$.
Now let us use formula \eqref{const-gener-eq} to define the pre-aisle $D^{\le 0}(X\times \P^r)$.
It follows from Lemma \ref{gener-lem} that 
$$D^{\le -1}_0(X\times \P^r)\sub D^{\le 0}(X\times \P^r)\sub D^{\le 0}_0(X\times \P^r).$$ 
By Theorem \ref{tilt-thm} this implies that $D^{\le 0}(X\times \P^r)$ extends to a $t$-structure.
Computing the right orthogonal to $D^{\le -1}(X\times \P^r)$ gives
\begin{equation}\label{proj-eq}
D^{\ge 0}(X\times\P^r)
=\{F\in D(X\times \P^r)\ |\ p_*(F\otimes \OO(n))\in D^{\ge 0}(X)\text{ for all }n\in\Z\}.
\end{equation}
It is easy to see that in this formula it is enough to consider
$n\gg 0$: one has to use the exactness of the Koszul complex 
\begin{equation}\label{Koszul-eq}
0\to\OO_{\P^r}(n-r)\to\sideset{}{^r}{\bigwedge} V\otimes\OO_{\P^r}(n-r+1)\to\ldots\to V\otimes\OO_{\P^r}(n)\to
\OO_{\P^r}(n+1)\to 0,
\end{equation}
where $V=H^0(\P^r,\OO_{\P^r}(1))$.
On the other hand, \eqref{const-gener-eq} implies 
that for every $F\in D^{\le 0}(X\times\P^r)$ one has $p_*(F(n))\in D^{\le 0}(X)$ for $n\gg 0$.
Hence, the left-hand side of \eqref{proj-eq2} is contained in its right-hand side.
Now the equality in \eqref{proj-eq2} follows from Lemma \ref{exact-fun-lem}(ii).

Since the constructed $t$-structure is invariant under tensoring with $\OO_{\P^r}(1)$,
by Theorem \ref{loc-thm} it is local over $\P^r$. 
\ed

In the situation of the above lemma let 
$\CC\sub D(X)$ (resp., $\CC_{\P^r}\sub D(X\times\P^r)$)
be the heart of the $t$-structure on $D(X)$ (resp., $D(X\times\P^r)$).
We have the following property (similar to Lemma 2.3.7 of \cite{AP}).

\begin{lem}\label{sur-lem} In the situation of Lemma \ref{proj-lem}
for every $F\in\CC_{\P^r}$ there exists a surjection $p^*G(n)\to F$ in $\CC_{\P^r}$
with $G\in\CC$ and $n\in\Z$.
\end{lem}

\Pf . Set $V=H^0(\P^r,\OO(1))$. First, we observe that for every $F\in\CC_{\P^r}$ the
natural map $V\otimes p_*(F(n))\to p_*(F(n+1))$ is surjective in $\CC$ for $n\gg 0$.
Indeed, this follows immediately from the exactness of the Koszul complex \eqref{Koszul-eq}.
Therefore, if we fix a sufficiently large $n$ then for $N>n$ the natural map
\begin{equation}\label{S-N-map}
S^{N-n}V\otimes p_*(F(n))\to p_*(F(N))
\end{equation}
is surjective. We claim that this implies surjectivity of the map $f: p^*p_*(F(n))(-n)\to F$ in $\CC_{\P^r}$.
Indeed, if $C=\coker(f)$ (the cokernel taken in $\CC_{\P^r}$)
then for $N\gg 0$ the object $p_*(C(N))\in\CC$ can be identified
with the cokernel of the map \eqref{S-N-map} in $\CC$. Hence, we get $p_*(C(N))=0$ for $N\gg 0$.
Therefore, $C=0$, i.e., $f$ is surjective.
\ed   
  
In the following lemma we compute the restriction of the above constant $t$-structure to $X\times\A^r$.
We denote by $(D^{\le 0}_{qc}(X), D^{\ge 0}_{qc}(X)$ the extension of our
$t$-structure from $D(X)$ to $D_{qc}(X)$, where $D^{\le 0}_{qc}(X)=\pra_{D_{qc}(X)}[[D^{\le 0}(X)]]$
and $D^{\ge 0}_{qc}(X)$ is the right orthogonal to $D^{\le -1}(X)$ in $D_{qc}(X)$
(see Lemma \ref{cocomplete-lem}).

\begin{lem}\label{aff-lem}
Let $(D^{\le 0}(X),D^{\ge 0}(X))$ be a close to Noetherian nondegenerate $t$-structure
on a scheme $X$. Then there exists a $t$-structure on $D(X\times\A^r)$, local over $\A^r$, given by
\begin{equation}\label{affine-eq}
D^{[a,b]}(X\times\A^r)=\{F\in D(X\times\A^r) \ |\ p_*(F)\in D^{[a,b]}_{qc}(X)\},
\end{equation}
where $p:X\times\A^r\to X$ is the projection. In addition we have 
\begin{equation}\label{affine-eq2}
D^{\le 0}(X\times\A^r)=\pra[p^*D^{\le 0}(X)].
\end{equation}
If the original $t$-structure on $D(X)$ is Noetherian then so is the constructed $t$-structure on $D(X\times\A^r)$.
\end{lem}

\Pf . Since the $t$-structure on $D(X\times\P^r)$ constructed in Lemma \ref{proj-lem} is local over
$\P^r$, it induces a $t$-structure on $D(X\times U)$ such that
for an open subset $U\sub\P^r$ the subcategory $D^{[a,b]}(X\times U)\sub D(X\times U)$ consists of restrictions
of objects of $D^{[a,b]}(X\times\P^r)$. This immediately gives \eqref{affine-eq2}.
Computing the orthogonal to $p^*D^{\le -1}(X)$ we get \eqref{affine-eq} for $[a,b]=[0,+\infty)$.
Also, using \eqref{affine-eq2} one easily checks that
$p_*D^{\le 0}(X\times\A^r)\sub D^{\le 0}_{qc}(X)$. Therefore, the functor
$p_*:D(X\times\A^r)\to D_{qc}(X)$ is $t$-exact with respect to our $t$-structures.
By Lemma \ref{exact-fun-lem}(i) this implies  
\eqref{affine-eq}. The assertion about Noetherian $t$-structures was proved in
Theorem 2.3.6 of \cite{AP}.
\ed

Below we will also use the following technical observation.

\begin{lem}\label{countable-lem}
Let $(D^{\le 0}(X),D^{\ge 0}(X))$ be a close to Noetherian nondegenerate $t$-structure
on a scheme $X$ of finite type over $k$. Extend it to a $t$-structure $D_{qc}(X)$ by setting
$D^{\le 0}_{qc}(X)=\pra_{D_{qc}(X)}[[D^{\le 0}(X)]]$. Then
for every $k$-vector space $V$ with a countable basis one has
$$V\otimes_k D^{[a,b]}(X)\sub D_{qc}^{[a,b]}(X).$$
\end{lem}

\Pf . It suffices to prove this for one vector space with an infinite countable basis.
Let us equip the categories $D(X\times\P^1)$ and $D(X\times\A^1)$ with constant $t$-structures
from Lemmas \ref{proj-lem} and \ref{aff-lem}. It is clear from the definition and the projection formula that the functor of pull-back $D(X)\to D(X\times\P^1)$ is $t$-exact. The restriction functor
$D(X\times\P^1)\to D(X\times\A^1)$ is also $t$-exact by locality of the constant $t$-structure
over $\P^1$. It follows 
that the pull-back functor $p^*:D(X)\to D(X\times\A^1)$ is $t$-exact, i.e.,
$p^*D^{[a,b]}(X)\sub D^{[a,b]}(X\times\A^1)$. Combining this with \eqref{affine-eq} we deduce
that $H^0(\A^1,\OO)\otimes_k D^{[a,b]}(X)\sub D_{qc}^{[a,b]}(X)$ as required.
\ed

Now we are ready to prove our main result about constant $t$-structures.

\begin{thm}\label{const-t-str-thm} 
Let $(D^{\le 0}(X),D^{\ge 0}(X))$ be a close to Noetherian nondegenerate $t$-structure
on a scheme $X$ of finite type over $k$. 

\noindent (i) 
For every scheme $S$ of finite type over $k$
we have a close to Noetherian nondegenerate $t$-structure on
$D(X\times S)$, local over $S$, such that
\begin{equation}\label{const-formula}
D^{[a,b]}(X\times S)=\{F\in D(X\times S)\ |\ p_*(F|_{X\times U})\in D^{[a,b]}_{qc}(X)
\text{ for every open affine }U\sub S\},
\end{equation}
where $p:X\times U\to X$ is the natural projection. 
If $S=\cup_i U_i$ is an open affine covering of $S$ then $F\in D^{[a,b]}(X\times S)$ iff 
$p_*(F|_{X\times U_i})\in D^{[a,b]}_{qc}(X)$ for every $i$.
If the original $t$-structure is Noetherian then so is the obtained $t$-structure on $D(X\times S)$.
The functor $p^*:D(X)\to D(X\times S)$ is $t$-exact with respect to the original $t$-structure
on $D(X)$ and the above $t$-structure on $D(X\times S)$.


\noindent (ii) Assume that $S$ is projective.
Then the above $t$-structure satisfies
\begin{equation}\label{proj-S-eq}
D^{[a,b]}(X\times S)=\{F\in D(X\times S)\ |\ p_*(F\otimes p_S^*L^n)\in D^{[a,b]}(X)\text{ for all }n\gg 0\},
\end{equation}
where $p_S:X\times S\to S$ is the projection, $L$ is an ample line bundle on $S$. 
\end{thm}

\Pf . (i) Assume first that $S$ is affine. 
Let us choose a closed embedding $S\hra\A^r$. Note that this embedding has finite Tor dimension
since $\A^r$ is smooth. Applying Theorem \ref{base-change-thm} to the constant
$t$-structure on $D(X\times\A^r)$ constructed in Lemma \ref{aff-lem}
we get an induced $t$-structure on $D(X\times S)$ with
$$D^{[a,b]}(X\times S)=\{F\in D(X\times S)\ |\ p_*(F)\in D^{[a,b]}_{qc}(X)\}.$$
Note that the right-hand side does not depend on a choice of an embedding $S\hra\A^r$.
Since $S$ is affine, from Theorem \ref{loc-thm} we see that this $t$-structure is local over $S$.
Let us extend this $t$-structure to a $t$-structure on $D_{qc}(X\times S)$
such that $D_{qc}^{\le 0}(X\times S)=\pra_{D_{qc}(X\times S)}[[D^{\le 0}(X\times S)]]$
(see Lemma \ref{cocomplete-lem}). We claim that 
\begin{equation}\label{const-qc-eq}
D_{qc}^{[a,b]}(X\times S)=\{F\in D_{qc}(X\times S)\ |\ p_*(F)\in D_{qc}^{[a,b]}(X)\}.
\end{equation}
First, we can check that the left-hand side is contained in the right-hand side. Indeed,
for $[a,b]=[-\infty,0]$ the right-hand-side of \eqref{const-qc-eq} is a cocomplete pre-aisle
containing $D^{\le 0}(X\times S)$. Hence, it also contains $D^{\le 0}_{qc}(X\times S)$.
On the other hand, we have the inclusion $p^*D^{\le -1}(X)\sub D^{\le -1}(X\times S)$
(since $p_*p^*F\simeq F\otimes H^0(S,\OO_S)$ and $D^{\le -1}_{qc}(X)$ is closed under
coproducts). Passing to right orthogonals in $D_{qc}(X\times S)$ we derive that
$p_*D_{qc}^{\ge 0}(X\times S)\sub D_{qc}^{\ge 0}(X)$.
Now the equality in \eqref{const-qc-eq} follows from Lemma \ref{exact-fun-lem}(i) applied to the functor
$p_*:D_{qc}(X\times S)\to D_{qc}(X)$.

Now let $j:U\hra S$ be an open affine subset. Then we have
a similar constant $t$-structure on $D(X\times U)$ and its extension to $D_{qc}(X\times U)$. 
From the above formulas it is clear that the functor 
$(\id_X\times j)_*:D_{qc}(X\times U)\to D_{qc}(X\times S)$ is $t$-exact. By adjunction we derive
that $(\id_X\times j)^*D_{qc}^{\le 0}(X\times S)\sub D_{qc}^{\le 0}(X\times U)$. On the other hand, 
$$D_{qc}^{\le 0}(X\times U)=(\id_X\times j)^*(\id_X\times j)_*D_{qc}^{\le 0}(X\times S)\sub 
(\id_X\times j)^*D_{qc}^{\le 0}(X\times S),$$
hence $D_{qc}^{\le 0}(X\times U)=(\id_X\times j)^*D_{qc}^{\le 0}(X\times S)$. It follows that
$$D_{qc}^{\le 0}(X\times U)=\pra_{D_{qc}(X\times U)}[[(\id_X\times j)^*D^{\le 0}(X\times S)]].$$
Since $(\id_X\times j)^*D^{\le 0}(X\times S)$ is an aisle in $D(X\times S)$ (by locality of
the constant $t$-structure over $S$), using Lemma \ref{cocomplete-lem} we get the following equality
of the aisles in $D(X\times U)$:
$$D^{\le 0}(X\times U)=(\id_X\times j)^*D^{\le 0}(X\times S).$$
Therefore, the corresponding $t$-structures coincide, so the functor
$(\id\times j)^*:D(X\times S)\to D(X\times U)$ is $t$-exact with respect to
the constant $t$-structures. This implies formula \eqref{const-formula} for affine $S$.

Note that if the original $t$-structure is Noetherian then so is
the constant $t$-structure on $D(X\times\A^r)$. By Theorem \ref{base-change-thm}
this also implies that the induced $t$-structure on $D(X\times S)$ is Noetherian for affine $S$.

Now let us consider the case of arbitrary $S$.
Let $S=\cup_i U_i$ be a finite open affine covering of $S$. 
By the preceding part of the proof, for every $i$ 
we have a $t$-structure on $D(X\times U_i)$ local over $U_i$. We claim that these $t$-structures agree on intersections $U_i\cap U_j$. Indeed, both the $t$-structures on $D(X\times (U_i\cap U_j))$, 
restricted from $X\times U_i$ and from $X\times U_j$, are local over $U_i\cap U_j$. Since
their restrictions to $X\times U$ agree for every open affine subset $U\sub U_i\cap U_j$, our claim
follows.
By Lemma \ref{local-glue-lem} the above $t$-structures on $D(X\times U_i)$ 
can be glued into a $t$-structure on $D(X\times S)$.
It is easy to see that it is still given by \eqref{const-formula}.

If the original $t$-structure is Noetherian then the constructed $t$-structure on $D(X\times S)$ will also be Noetherian (this immediately reduces to the affine case considered above).
Since our construction preserves inclusions between pre-aisles, 
we derive that the constant $t$-structure on $D(X\times S)$ is close to Noetherian.

Finally, let us check that the functor $p^*:D(X)\to D(X\times S)$ is $t$-exact. If $S$ is affine
then we have $p_*p^*F\simeq H^0(S,\OO_S)\otimes F$, so the assertion follows from
Lemma \ref{countable-lem}. The general case is deduced easily 
by covering $S$ with
open affine subsets and using locality of our $t$-structure on $D(X\times S)$ over $S$.

\noindent
(ii) Since the constant $t$-structure on $D(X\times S)$ is local over $S$, it is invariant under
tensoring with the pull-back of any line bundle on $S$. Hence, for any integer $d>0$ and
$F\in D(X\times S)$ we have
$F\in D^{[a,b]}(X\times S)$ iff $F\otimes p_S^*L^i\in D^{[a,b]}(X\times S)$ for $i=0,\ldots,d-1$.
Therefore, it is enough to prove \eqref{proj-S-eq} for $L^d$
instead of $L$, so we can assume that $L$ is very ample.
Let $i:S\to\P^r$ be a closed embedding such that $L=i^*\OO(1)$.
Applying Theorem \ref{base-change-thm} to the constant $t$-structure on $D(X\times\P^n)$
we derive that the right-hand side of \eqref{proj-S-eq} 
gives a $t$-structure on $D(X\times S)$ (automatically local over $S$ by Theorem \ref{loc-thm}).
Considering the standard open covering of $\P^n$ by the affine pieces,
it is easy to see that the restrictions of this $t$-structure on $D(X\times S)$ to the induced open affine pieces of $S$ agree with the $t$-structures constructed in (i). Hence,
it coincides with the $t$-structure given by \eqref{const-formula}.
\ed

In the case of an affine base $S=\Spec(A)$ the heart of the constant $t$-structure on
$D(X\times S)$ has a natural description in terms of $A$-modules in the heart of the $t$-structure on $D_{qc}(X)$. More precisely, let $\CC\sub D(X)$ be the heart of the original $t$-structure on $D(X)$, 
and let $\CC_{qc}\sub D_{qc}(X)$ be the heart of the corresponding $t$-structure on $D_{qc}(X)$
(such that $D_{qc}^{\le 0}(X)=\pra_{D_{qc}(X)}[[D^{\le 0}(X)]]$).
Recall that an $A$-{\it module} in $\CC_{qc}$ is an object $F\in\CC_{qc}$ equipped with
a homomorphism of algebras $A\to\End(F)$ (see \cite{Po}). 
They form an abelian category that we will denote
by $A-\mod-\CC_{qc}$. Note that since $A$ has a countable basis as a $k$-vector space,
by Lemma \ref{countable-lem} for every $F\in\CC$ we have $F\otimes_k A\in\CC_{qc}$. 
Let us say that an $A$-module $M$ in $\CC_{qc}$ is {\it finitely presented} if 
there exists a pair of objects $F_0,F_1\in\CC$ and a morphism of free $A$-modules
$f:F_1\otimes_k A\to F_0\otimes_k A$
in $A-\mod-\CC_{qc}$ such that $M\simeq\coker(f)$. 

\begin{prop}\label{heart-prop} 
Keep the assumptions of Theorem \ref{const-t-str-thm}. Assume in addition
that $S=\Spec(A)$ for a finitely generated $k$-algebra $A$. Then the heart $\CC_S$
of the constant $t$-structure on $D(X\times S)$ 
is equivalent to the category of finitely presented $A$-modules in $\CC_{qc}$.
\end{prop}

\Pf . We have the exact functor $p_*:\CC_S\to \CC_{qc}$. 
Furthermore, for every $F\in\CC_S$ the object $p_*F\in\CC_{qc}$ has a natural 
$A$-module structure given by the homomorphism 
$$A\to H^0(X\times S,\OO)\to\End(F)\to\End(p_*F).$$
Thus, we obtain an exact functor $p_*:\CC_S\to A-\mod-\CC_{qc}$.
It sends $p^*G$, where $G\in\CC$, to
the free $A$-module $p_*p^*G\simeq G\otimes_k A$. We claim that for every $F\in\CC_S$ there
exists a surjection $p^*G\to F$ in $\CC_S$ with $G\in\CC$. Indeed, assume first that $S=\A^r$.
Then we can extend every $F\in\CC_{\A^r}$ to an object $\ov{F}\in\CC_{\P^r}$ such that
$\ov{F}|_{X\times\A^r}\simeq F$. 
Applying Lemma \ref{sur-lem} to $\ov{F}$ and restricting the obtained surjection to $X\times\A^r$
we deduce our claim for $S=\A^r$. For general $S$ let us consider a closed embedding
$i:S\to\A^r$. Then there exists a surjection of the required type for $(\id_X\times i)_*F$ in $\CC_{\A^r}$.
It remains to restrict it to $S$ and use Lemma \ref{closed-emb-lem} together with the fact that
$(\id_X\times i)^*$ is right $t$-exact (as the left adjoint of the $t$-exact functor $(\id_X\times i)_*$). 
This proves our claim. It follows that every object in $\CC_S$
can be represented as the cokernel of a morphism $p^*G_1\to p^*G_0$, where $G_0,G_1\in\CC$.
Hence, the $A$-module $p_*F$ is finitely presented for every $F\in\CC$.
Next, we have a natural isomorphism 
$$\Hom_{\CC_S}(p^*G,F)\simeq\Hom_{\CC_{qc}}(G,p_*F)\simeq
\Hom_{A-\mod-\CC_{qc}}(G\otimes_k A,p_*F)$$
for $G\in\CC$, $F\in\CC_S$. Representing arbitrary $F\in\CC_S$ as the cokernel of a morphism $p^*G_1\to p^*G_0$
with $G_0,G_1\in\CC$ we deduce that 
$$\Hom_{\CC_S}(F,F')\simeq\Hom_{A-\mod-\CC_{qc}}(p_*F,p_*F')$$
for $F,F'\in\CC_S$. It is also clear that every finitely presented $A$-module in $\CC_{qc}$
is in the essential image of the functor $p_*$, so our assertion follows.
\ed

\subsection{Localization}

Let $f:X\to S$ be a morphism, where $S$ is quasiprojective. 
As an application of our techniques, we show that some $t$-structures on $D(X)$ 
naturally give rise to new $t$-structures on $D(X)$, local over $S$.

\begin{thm}\label{localization-thm} 
Let $L$ be an ample line bundle on $S$, and let
$(D^{\le 0}(X),D^{\ge 0}(X))$ be a nondegenerate close to Noetherian $t$-structure
on $D(X)$ such that tensoring with $f^*L$ is right $t$-exact, i.e., $f^*L\otimes D^{\le 0}(X)\sub D^{\le 0}(X)$. Then there exists a $t$-structure on $D(X)$, local over $S$, given by
\begin{equation}\label{D-f-eq}
D_f^{[a,b]}(X)=\{ F\in D(X)\ |\ F\otimes f^*L^n\in D^{[a,b]}(X)\text{ for all }n\gg 0\}.
\end{equation}
We also have
\begin{equation}\label{D-f-eq2}
D_f^{[a,b]}(X)=\{ F\ |\ j_*j^*F\in D^{[a,b]}_{qc}(X)\text{ for every open }U\sub S,
\text{ where }j:f^{-1}(U)\hra X\}.
\end{equation}
If the original structure is Noetherian then so is the new one.
\end{thm}

\Pf . First, let us check that \eqref{D-f-eq} is a $t$-structure on $D(X)$, local over $S$.
Let us denote the right-hand side of \eqref{D-f-eq} by $D_{f,L}^{[a,b]}(X)$.
We claim that it is enough to prove that $D_{f,L^d}^{[a,b]}(X)$ is a $t$-structure for
some $d>0$. Indeed, by Theorem \ref{loc-thm} this $t$-structure is local over $S$.
Hence, it is stable under tensoring with $f^*L$, so we have $F\in D_{f,L^d}^{[a,b]}(X)$
iff $F\otimes f^*L^i\in D_{f,L^d}^{[a,b]}(X)$ for $i=0,\ldots,d-1$. Therefore,
$D_{f,L}^{[a,b]}(X)=D_{f,L^d}^{[a,b]}(X)$, and our claim follows.
Thus, we can assume that $L$ is very ample.
Let $\iota:S\to\P^r$ be the locally closed embedding such that $\iota^*\OO_{\P^r}(1)=L$.
Since the right-hand-side of \eqref{D-f-eq} depends only on $f^*L$, it is enough to
prove that \eqref{D-f-eq} gives a $t$-structure after replacing the data 
$(S,f,L)$ with $(\P^r,\iota\circ f,\OO_{\P^r}(1))$. Thus, we can assume
from the beginning that $S=\P^r$ and $L=\OO_{\P^r}(1)$.
Consider the closed embedding $i=(\id,f):X\to X\times \P^r$. 
We claim that under our assumptions the functor of tensoring with $i_*\OO_X$ on $D(X\times \P^r)$
is right $t$-exact with respect to the constant $t$-structure.
Indeed, by \eqref{const-gener-eq} it suffices to prove the inclusions
$$i_*\OO_X\otimes (D^{\le 0}(X)\boxtimes\OO_{\P^r}(n_0))\sub D^{\le 0}(X\times\P^r)$$
for all $n_0\in\Z$. For $F\in D^{\le 0}(X)$ we have
$$i_*\OO_X\otimes (F\boxtimes\OO_{\P^r}(n_0))\simeq i_*(F\otimes f^*\OO_{\P^r}(n_0)).$$
Tensoring this with $\OO_{\P^r}(n)$ and pushing forward to $X$ gives
$F\otimes f^*\OO_{\P^r}(n_0+n)$ which belongs to $D^{\le 0}(X)$ for $n\gg 0$ by the assumption.
Note also that $i$ has finite Tor dimension since $\P^r$ is smooth.
Thus, we can apply Proposition \ref{finite-prop}. The 
induced $t$-structure on $D(X)$ will be given by \eqref{D-f-eq}.

From \eqref{D-f-eq} we deduce that the push-forward with respect to the closed embedding
$(\id,f):X\to X\times S$ is $t$-exact with respect to our new $t$-structure on $D(X)$ and
the constant $t$-structure $(D^{\le 0}(X\times S),D^{\ge 0}(X\times S))$
on $D(X\times S)$ induced by the old $t$-structure on $D(X)$.
Hence, 
$$D^{[a,b]}_f(X)=\{ F\in D(X)\ |\ (\id,f)_*F\in D^{[a,b]}(X\times S)\}.$$
Now \eqref{D-f-eq2} follows from \eqref{const-formula}.
\ed

\begin{rem} The assumption $f^*L\otimes D^{\le 0}(X)\sub D^{\le 0}(X)$ in the above
theorem is equivalent (by passing to right orthogonals) to the condition
$f^*L^{-1}\otimes D^{\ge 0}(X)\sub D^{\ge 0}(X)$. Therefore, the formula for
$D_f^{\ge 0}(X)$ can be rewritten as
$$D_f^{\ge 0}(X)=\{ F\in D(X)\ |\ F\otimes f^*L^n\in D^{\ge 0}(X)\text{ for all }n\in\Z\},$$
so that $D_f^{\ge 0}(X)=\cap_{n\in\Z} D_n^{\ge 0}(X)$, where
$D_n^{\ge 0}(X)=f^*L^{-n}\otimes D^{\ge 0}(X)$. Note that we have a chain of inclusions 
$$\ldots\supset D_n^{\ge 0}\supset D_{n+1}^{\ge 0}\supset\ldots$$
Thus, Theorem \ref{localization-thm} can be viewed as an example of a ``limiting" $t$-structure,
like Corollary \ref{lim-cor}. Note also that if we apply this theorem to the glued $t$-structures 
$(D_n^{\le 0}(X\times\P^r), D_n^{\ge 0}(X\times\P^r))$ on
$D(X\times\P^r)$ (associated with a $t$-structure on $D(X)$, see section \ref{const-t-str-sec})
and take $L=\OO_{\P^r}(1)$ then it will produce the constant $t$-structure
on $D(X\times\P^r)$.
\end{rem}

In the case when $X$ is quasiprojective and
$f$ is the identity morphism the above construction produces a local $t$-structure on $D(X)$. 

\begin{cor} Assume that $X$ is quasiprojective. Let $L$ be an ample line bundle on $X$.
Let also $(D^{\le 0}(X),D^{\ge 0}(X))$ be a nondegenerate close to Noetherian $t$-structure
on $D(X)$ such that $L\otimes D^{\le 0}(X)\sub D^{\le 0}(X)$.
Then for every smooth point $x\in X$ the structure sheaf $\OO_x$ has only one
nonzero cohomology object with respect to $(D^{\le 0}(X),D^{\ge 0}(X))$.
\end{cor}

\Pf . By Theorem \ref{base-change-thm} the corresponding local $t$-structure
$(D_{\id}^{\le 0},D_{\id}^{\ge 0})$ is compatible with some $t$-structure on $D(x)$ for every
smooth point $x\in X$ (so that
the push-forward functor with respect to the embedding $x\hra X$ is $t$-exact). 
Hence, $\OO_x$ has only one nonzero cohomology with respect to 
$(D_{\id}^{\le 0},D_{\id}^{\ge 0})$. It remains to use the fact that $L^n\otimes\OO_x\simeq\OO_x$.
\ed 

\begin{rem} The condition that $L\otimes D^{\le 0}(X)\sub D^{\le 0}(X)$ is crucial
in the above corollary. Without it the assertion may be wrong even
for a nondegenerate Noetherian $t$-structure. For example, let $X$ be a K3 surface,
$C\sub X$ a $(-2)$-curve. Take $(D^{\le 0}(X),D^{\ge 0}(X))$ to be the image of the standard
$t$-structure under the reflection functor $T_{\OO_C}^{-1}$ (see \cite{ST}). Then for $x\in C$ the
structure sheaf $\OO_x$ will have two nontrivial cohomology objects.
\end{rem}

\subsection{Invariance under a connected group of autoequivalences}\label{invar-sec}

In this section we assume that our ground field $k$ is algebraically closed.

Recall that if $K\in D(X\times Y)$ is an object of finite Tor-dimension, such that
its support is proper over $Y$, then it induces an exact functor
$$\Phi_K:D(X)\to D(Y): F\mapsto p_{2*}(p_1^*F\otimes K),$$
where $p_i$ are the projections from $X\times Y$ to its factors.
We say that $K$ is the {\it kernel} giving the functor $\Phi_K$.
It follows from the theorem of Orlov in \cite{Orlov} that if $X$ is a smooth projective variety then every exact autoequivalence of $D(X)$ is given by some kernel. 

Let us denote by $\Autoeq D(X)$ the group of (isomorphism classes of) exact autoequivalences of
$D(X)$. By an action of a group $G$ on $D(X)$ we mean a homomorphism
$G\to \Autoeq D(X): g\mapsto \Phi_g$.
In the case when $G$ is an algebraic group there is a natural way to strengthen this definition
by requiring the existence of a family of kernels.

\begin{defi}
We say that an algebraic group $G$ acts on $D(X)$ by {\it kernel autoequivalences} if
we are given a homomorphism $G\to\Autoeq D(X):g\mapsto \Phi_g$, and
an object $K\in D(G\times X\times X)$ of finite Tor dimension with the support
proper over $G\times X$ (with respect to the projection $p_{13}$), such that 
for every $g\in G(k)$ we have $\Phi_g=\Phi_{K_g}$, where $K_g=K|_{g\times X\times X}$.
\end{defi}

For example, the Poincar\'e line bundle $\PP$ on $\Pic^0(X)\times X$ gives rise to an action
of $\Pic^0(X)$ on $D(X)$ by kernel autoequivalences. Namely, we should take
$K=(\id\times\De)_*\PP\in D(\Pic^0(X)\times X\times X)$, where $\De:X\to X\times X$ is
the diagonal.

\begin{thm}\label{invar-thm} 
Let $(D^{\le 0}(X),D^{\ge 0}(X))$ be a nondegenerate Noetherian $t$-structure 
on $D(X)$. Assume that a connected smooth algebraic group $G$ acts on $D(X)$ by
kernel autoequivalences. Then $(D^{\le 0}(X),D^{\ge 0}(X))$ is invariant under this action.
\end{thm}

\Pf . Let $\CC=D^{\le 0}(X)\cap D^{\ge 0}(X)$. For $F\in\CC$ consider the object
$K*F\in D(G\times X)$ defined by
$$K*F=p_{13*}(K\otimes p_2^*F),$$
where $p_{ij}$ and 
$p_i$ are projections from $G\times X\times X$, $K$ is the kernel defining the action of
$G$. Then $(K*F)|_{g\times X}\simeq \Phi_g(F)$.
In particular, $(K*F)|_{e\times X}\simeq F\in\CC$, where $e\in G$ is the neutral element. By the open heart property (see Proposition \ref{open-heart-prop}) this implies
that there exists an open neighborhood $U$ of $e$ in $G$ such that
$(K*F)|_{U\times X}$ belongs to the heart of the constant $t$-structure on $D(U\times X)$.
Since for any $g\in U$ the restriction functor $D(U\times X)\to D(\{g\}\times X)$ is right $t$-exact
(as the left adjoint to the $t$-exact push-forward functor),
this implies that $\Phi_g(F)\in D^{\le 0}(X)$ for all $g\in U$. Thus, the functors
$\Phi_g$ are right $t$-exact for all $g\in U$. It follows that for $g\in U\cap U^{-1}$
the functors $\Phi_g$ are $t$-exact.
Hence, the set of $g$ such that $\Phi_g$ is $t$-exact is an open subgroup in $G$,
so it is equal to $G$.
\ed

\begin{cor}\label{invar-cor} 
Assume $X$ is smooth and projective.
Let $\Sigma$ be a connected component in the space of numerical stability conditions
on $D(X)$ such that the corresponding subspace $V(\Sigma)\sub(\NN(X)\otimes\C)^*$
is defined over $\Q$ (see \cite{Bridge1}). Then any stability in $\Sigma$ is invariant under the action
of a connected group of kernel autoequivalences.
\end{cor}

\Pf . Indeed, for a stability with Noetherian $\PP(0,1]$ this follows from Theorem \ref{invar-thm}.
Since the set of such stabilities is dense in $\Sigma$
and autoequivalences act by isometries, the general case follows.
\ed

\end{document}